\documentclass{amsart}%
\usepackage{amsfonts}
\usepackage{amsmath}
\usepackage{amssymb}
\usepackage{graphicx}
\usepackage{xcolor,ulem}
\usepackage[colorlinks, citecolor = blue, pagebackref]{hyperref}%
\providecommand{\U}[1]{\protect\rule{.1in}{.1in}}
\newtheorem{theorem}{Theorem}
\theoremstyle{plain}

\newtheorem{corollary}{Corollary}

\newtheorem{example}{Example}

\newtheorem{lemma}{Lemma}

\newtheorem{proposition}{Proposition}
\newtheorem{remark}{Remark}

\numberwithin{equation}{section}
\begin{document}
\title[Majorization revisited: Comparison of norms in interpolation scales]{Majorization revisited: Comparison of norms in interpolation scales\\}
\author{Sergey V. Astashkin}
\address[S.V. Astashkin]{Samara National Research University, Moskovskoye shosse 34,
443086, Sama\-ra, Russia}
\email{astash56@mail.ru}
\urladdr{https://ssau.ru/english/staff/334603-astashkin-sergey-v}
\author{Konstantin V. Lykov}
\address[K.V. Lykov]{Image Processing Systems Institute --- Branch of the Federal
Scientific Researh Centre "Crystallography and Photonics" of Russian Academy
of Sciences, Molodogvardejskaya st. 151, 443001 Samara, Russia and Faculty of
Information Technology, Samara National Research University, 443086 Samara, Russia}
\email{alkv@list.ru}
\author{Mario Milman}
\address[M. Milman]{Instituto Argentino de Matem\'{a}tica\\
Buenos Aires, Argentina}
\email{mario.milman@icloud.com}
\urladdr{https://sites.google.com/site/mariomilman/}
\thanks{The authors collaboration for this research was initially supported by a grant
from the Simons Foundation ($\backslash$\#207929 to Mario Milman). The work of
the first named author was completed as a part of the implementation of the
development program of the Volga Region Scientific and Educational
Mathematical Center (agreement no. 075-02-2021-1393). The second named author
has been supported by the RFBR grant 19-29-09054.}
\date{}

\begin{abstract}
We reformulate, modify and extend a comparison criteria of $L^{p}$ norms
obtained by Nazarov-Podkorytov and place it in the general setting of
interpolation theory and majorization theory. In particular, we give norm
comparison criteria for general scales of interpolation spaces, including
non-commutative $L^{p}$ and Lorentz spaces. As an application, we extend the
classical Ball's integral inequality, which lies at the basis of his famous
result on sections of the $n-$dimensional unit cube.

\end{abstract}
\maketitle
\tableofcontents

\section{Introduction: Background and plan of the paper}

The starting point of our research is the integral inequality (cf. \cite{B1},
\cite{B2}, \cite{np})%
\begin{equation}
\int_{-\infty}^{+\infty}\left\vert \frac{\sin\pi x}{\pi x}\right\vert
^{p}\,dx<\sqrt{\frac{2}{p}}\;\;\mbox{for every}\;\;p>2, \label{intro1}%
\end{equation}
(for $p=2$, it turns into an identity), which lies at the basis of Ball's
famous result on sections of the $n-$dimensional unit cube (cf. \cite{B1}, for
more background on these problems we refer to the recent survey \cite{klar}).
The original proof of (\ref{intro1}) involved a large number of unappealing
numerical computations. Subsequently, Nazarov-Podkorytov \cite{np} gave a new
interesting proof introducing a variant of the classical method of comparison
of $L^{p}$-norms of functions through pointwise comparison of their
distribution functions. More recently, this approach was intensively developed
and used to establish sharp Khintchine-type inequalities (cf. \cite{BC},
\cite{Ko}, \cite{Mordhorst}, \cite{CKT-20} and the references therein) as well
as sharp comparison of moments for various classes of random variables
appearing in a geometric context \cite{Eske}. This was an additional
motivation for us to work on reformulating the Nazarov-Podkorytov approach in
the more general context of norm comparisons of elements belonging to an
interpolation scale. In particular, we introduce new majorization criteria
connected with interpolation that extends and clarifies the methods of
\cite{np}. As an application, we obtain general variants of Ball's $L^{p}%
$-integral inequality \eqref{intro1}, and improvements on related recently
published inequalities (cf. \cite{CKT-20} and the references therein).

The method of \cite{np} starts with a reformulation of inequality
(\ref{intro1}). Using the Euler-Poisson integral formula, $\int_{-\infty
}^{+\infty}e^{-\pi x^{2}}\,dx=1$, \ Ball's inequality (\ref{intro1}) can be
rewritten as follows. Let $\;$%
\begin{equation}
f(x):=e^{-\pi x^{2}/2},\;\;g(x):=\left\vert \frac{\sin\pi x}{\pi x}\right\vert
,\;\;\;x\in(0,\infty), \label{f&g}%
\end{equation}
then we want to prove that for $p>2$
\begin{equation}
\int_{0}^{+\infty}g(x)^{p}\,dx<\int_{0}^{+\infty}f(x)^{p}\,dx, \label{tres}%
\end{equation}
with equality when $p=2.$ At this point our problem is to compare the $L^{p}%
$-norms of the given functions $f$ and $g$. We shall now review some of the
considerations involved when comparing $L^{p}$-norms.

Let $(\Omega,\mu)$ be a measure space\footnote{We assume our measure spaces
are \textquotedblleft resonant\textquotedblright\ in the sense of \cite{BS}.}.
Given a measurable function $f:\Omega\rightarrow\mathbb{R}$, its distribution
function is defined by
\[
\lambda_{f}(\tau)=\mu\{\left\vert f\right\vert >\tau\},\;\;\tau>0.
\]
Suppose that there exists $C>0,$ such that $f,g$ satisfy
\begin{equation}
\lambda_{g}(\tau)\leq C\lambda_{f}(\tau),\text{ for all }\tau>0, \label{unou}%
\end{equation}
then we can easily conclude that
\[
\left\Vert g\right\Vert _{p}\leq C^{1/p}\left\Vert g\right\Vert _{p},\text{
for all }p\in\lbrack1,\infty)
\]
(in what follows, $\Vert\cdot\Vert_{p}$ denotes the norm of $L^{p}(\Omega)$).

It turns out, however, that in order to successfully obtain (\ref{tres}) it is
not possible to apply this elementary comparison criteria to the functions $f$
and $g$ that were defined in (\ref{f&g}). Indeed, on the one hand, we would
need for the constant $C$ that appears in (\ref{unou}) to be equal to $1,$
while on the other, since we know that $\left\Vert f\right\Vert _{2}%
^{2}-\left\Vert g\right\Vert _{2}^{2}=0,$ it follows that%
\[
\int_{0}^{\infty}2\tau[\lambda_{f}(\tau)-\lambda_{g}(\tau)]\,d\tau=0.
\]
Therefore, the fact that the function $\lambda_{f}(\tau)-\lambda_{g}(\tau)$,
$\tau>0$, is not identically zero implies the existence of $\tau_{0}>0,$ where
this function changes signs and, consequently, (\ref{unou}) cannot hold (with
$C=1$) for all $\tau>0$.

But not all is lost, in \cite{np} the authors show that for the functions
defined in (\ref{f&g}), the change of signs is unique, and propose a modified
general comparison criteria, replacing the assumption (\ref{unou}) by a weaker
condition to accommodate changes of signs. Let $\ f$ and $g$ be two arbitrary
measurable functions and, moreover, suppose that there exists $0<\tau
_{0}<\infty$ such that
\begin{equation}
\lambda_{f}(\tau)\left\{
\begin{array}
[c]{cc}%
\leq & \lambda_{g}(\tau)\text{ for }\tau<\tau_{0}\\
\geq & \lambda_{g}(\tau)\text{ for }\tau>\tau_{0}.
\end{array}
\right.  \label{distrib}%
\end{equation}
Furthermore, suppose that there exists $p_{0}\geq1$ such that
\begin{equation}
\int_{\Omega}(\left\vert f\right\vert ^{p_{0}}-\left\vert g\right\vert
^{p_{0}})d\mu\geq0. \label{penso}%
\end{equation}
Under these assumptions,\ a version of the Nazarov-Podkorytov Lemma
\cite[p.~250-251]{np} asserts that for all $p>p_{0},$ such that $|f|^{p}%
-|g|^{p}\in L^{1}(\Omega)$, we have
\begin{equation}
\int_{\Omega}(\left\vert f\right\vert ^{p}-\left\vert g\right\vert ^{p}%
)d\mu\geq0. \label{penso2}%
\end{equation}
It is shown in \cite{np} that, indeed, for the functions $f$ and $g,$ defined
in (\ref{f&g}), there exists a unique $\tau_{0}>0$ such that (\ref{distrib})
holds. Consequently, since (\ref{penso}) is valid (with equality) for
$p_{0}=2$, the above result allows us to obtain (\ref{tres}), from which
(\ref{intro1}) follows.

The proof given in \cite{np} that (\ref{distrib}) and (\ref{penso}) imply
(\ref{penso2}), is both elementary and ingenious. As it turns out, two ideas
are crucial: the use of distribution functions and the formulation of the
result in terms of differences. It will be useful for our future development
to go over the details.

Suppose that $f$ and $g$ are two arbitrary functions that satisfy the
assumptions (\ref{distrib}) and (\ref{penso}); furthermore, suppose that
$p>p_{0}$ is such that $|f|^{p}-|g|^{p}\in L^{1}(\Omega).$ Then, by
computation, we have
\begin{align*}
\int_{\Omega}\left(  \left\vert f\right\vert ^{p}-\left\vert g\right\vert
^{p}\right)  \,d\mu &  =p\int_{0}^{\infty}\left(  \lambda_{f}(s)-\lambda
_{g}(s)\right)  s^{p-1}\,ds\\
&  =p\tau_{0}^{p-1}\int_{0}^{\infty}\left(  \lambda_{f}(s)-\lambda
_{g}(s)\right)  \left(  \frac{s}{\tau_{0}}\right)  ^{p-1}\,ds.
\end{align*}
Therefore \eqref{penso2} will be proved once we verify that%
\[
\phi_{f,g}(p):=\int_{0}^{\infty}\left(  \lambda_{f}(s)-\lambda_{g}(s)\right)
\left(  \frac{s}{\tau_{0}}\right)  ^{p-1}ds\geq0.
\]
To accomplish this goal, we compare $\phi_{f,g}(p)$ with the known nonnegative
quantity $\phi_{f,g}(p_{0})$ (cf. (\ref{penso})), and write%
\begin{equation}
\phi_{f,g}(p)-\phi_{f,g}(p_{0})=\int_{0}^{\infty}\left(  \lambda
_{f}(s)-\lambda_{g}(s)\right)  \left(  \left(  \frac{s}{\tau_{0}}\right)
^{p-1}-\left(  \frac{s}{\tau_{0}}\right)  ^{p_{0}-1}\right)  ds.
\label{penso1}%
\end{equation}
To analyze the sign of the integrand we split $(0,\infty)=(0,\tau_{0}%
)\cup(\tau_{0},\infty).$ Since $\ln(\frac{s}{\tau_{0}})>0$ iff $s>\tau_{0},$
it holds that $\gamma(p):=\left(  \frac{s}{\tau_{0}}\right)  ^{p-1}%
=e^{(p-1)\ln(\frac{s}{\tau_{0}})},$ as a function of $p\geq1,$ is increasing
for each $s\in(\tau_{0},\infty)$ and is decreasing for $s\in(0,\tau_{0}).$
This fact, combined with (\ref{distrib}), shows that the factors $\left(
\lambda_{f}(s)-\lambda_{g}(s)\right)  $ and $\left(  \left(  \frac{s}{\tau
_{0}}\right)  ^{p-1}-\left(  \frac{s}{\tau_{0}}\right)  ^{p_{0}-1}\right)  ,$
have the same signs in each of the two intervals under consideration.
Therefore\footnote{An alternative method to prove that $\phi_{f,g}(p)$ is an
increasing function of $p$ is to show that for $p>p_{0}$
\[
\frac{d}{dp}(\phi_{f,g}(p))=\int_{0}^{\infty}\left(  \lambda_{f}%
(s)-\lambda_{g}(s)\right)  \left(  \frac{s}{\tau_{0}}\right)  ^{p-1}\ln
(\frac{s}{\tau_{0}})ds\geq0.
\]
In fact, the analysis of the signs of the factors $\left(  \lambda
_{f}(s)-\lambda_{g}(s)\right)  $ and $\ln(\frac{s}{\tau_{0}})$ is exactly the
same as the one provided above.},%
\[
\phi_{f,g}(p)-\phi_{f,g}(p_{0})\geq0.
\]
Finally, since $\phi_{f,g}(p_{0})\geq0,$ we conclude that
\begin{equation}
\phi_{f,g}(p)\geq0,\text{ for }p>p_{0}, \label{cuatro}%
\end{equation}
which is equivalent to \eqref{penso2}.

In the first part of the paper we extend the comparison method of $L^{p}%
$-norms of \cite{np} to the context of $L(p,q)$-spaces and non-commutative
Schatten ideals. As an application we obtain a new $L(p,q)$ version of the
Ball's inequality (\ref{intro1}). To describe the result we recall that for
$0<p,q<\infty,$ the Lorentz $L(p,q)$-spaces can be defined using the
quasi-norms%
\begin{equation}
\Vert f\Vert_{L(p,q)}:=\Big(\int_{0}^{\infty}(\lambda_{f}(s))^{q/p}%
\,d(s^{q})\Big)^{1/q}. \label{quasi-norm}%
\end{equation}
Equivalently, by a change of variables argument, we can write
(\ref{quasi-norm}) in terms of non-increasing rearrangements\footnote{The
non-increasing rearrangement $f^{\ast}$ of a measurable function
$f:\,\Omega\rightarrow\mathbb{R}$ {can be defined by the formula%
\[
f^{\ast}(s):=\inf\{\tau>0:\,\lambda_{f}(\tau)\leq s\},\text{ \ }s\in
(0,\mu(\Omega)).
\]
}}:
\begin{equation}
\Vert f\Vert_{L(p,q)}=\Big(\int_{0}^{\infty}f^{\ast}(s)^{q}\,d(s^{q/p}%
)\Big)^{1/q}. \label{quasi-norm*}%
\end{equation}

In Section \ref{s3} we obtain necessary and sufficient conditions on the
parameters $p$ and $q$ under which a version of the Nazarov-Podkorytov Lemma
holds for the $L(p,q)$-spaces. More precisely, let $0<{p_{0}},{q_{0}}<\infty$
be fixed. Let us denote by ${\mathcal{A}}({p_{0}},{q_{0}})$ the set of all
pairs $({p},{q})$, $0<{p},{q}<\infty$, such that for arbitrary measurable
functions $f,g$, that satisfy condition \eqref{distrib} for some $\tau_{0}>0$
and, moreover,
\[
\int_{0}^{\infty}(f^{\ast}(t)^{{q_{0}}}-g^{\ast}(t)^{{q_{0}}})\,d(t^{{q_{0}%
}/{p_{0}}})\geq0,
\]
from
\[
\ (f^{\ast}(t)^{{q}}-g^{\ast}(t)^{{q}})t^{{q}/{p}-1}\in L_{1}(0,\infty)
\]
it follows
\[
\int_{0}^{\infty}(f^{\ast}(t)^{{q}}-g^{\ast}(t)^{{q}})\,d(t^{{q}/{p}})\geq0.
\]
In Theorem \ref{prop of Dlta}, we show that
\[
{\mathcal{A}}({p_{0}},{q_{0}})=\{({p},{q}):\,{p}\geq{p_{0}},{q}\geq{q_{0}}%
,{q}/{p}\leq{q_{0}}/{p_{0}}\}.
\]

As an application, in Theorem \ref{teor<} of Section \ref{sec:integralineq} we
prove the following version of Ball's inequality (\ref{intro1}): if $q\geq2$
and $0<\alpha\leq1,$ then
\begin{equation}
\int_{0}^{+\infty}\left(  \left(  \frac{\sin t}{t}\right)  ^{\ast}\right)
^{{q}}t^{\alpha-1}\,dt\leq\frac{1}{2}q^{-\alpha/2}(2\pi)^{\alpha/2}%
\Gamma(\alpha/2), \label{show1}%
\end{equation}
where $\Gamma(t)$ is the $\Gamma$-function. We also compare (\ref{show1}) with
recent related inequalities obtained by Chasapis, K\"{o}nig and Tkocz in
\cite{CKT-20}.

One can formulate an analogous comparison method also for non-commutative
$L^{p}$-spaces. Given two compact operators, $F,G,$ say, acting on a Hilbert
space, we compare their Schatten $L^{p}$-norms by comparing their
non-commutative distribution functions or rearrangements (cf. \cite{Grot},
\cite{pee}) and we can easily state comparison criteria of \ \textquotedblleft
Nazarov-Podkorytov type\textquotedblright\ to prove that a suitable version of
(\ref{penso2}) holds. This is done in Section \ref{s4}.

In order to reach the widest possible audience we have chosen to present our
results concerning $L(p,q)$-spaces and Schatten ideals without explicit use of
the notation of interpolation theory. But simply remark here that since
through the use of the functionals of real interpolation one can recover
distribution functions and rearrangements, it is not difficult to see that
these results can be unified in the general framework of scales of
interpolation spaces. This point of view is developed in the next sections.

Interpolation theory plays a more substantial r\^{o}le in our development in a
different direction. In fact, the connection between interpolation and
majorization led us to discover that the crossing conditions of distribution
functions of \cite{np} imply majorization conditions which, in turn, via the
Calder\'{o}n-Mityagin interpolation theorem, allow us to prove much stronger
results and inequalities. Indeed, the results in the applications are
apparently out of the reach of previous methods.

In Section \ref{K-method} we extend (\ref{intro1}) to the class of
rearrangement invariant spaces (we use their definition given in \cite{LT}) as follows.

\begin{theorem}
\label{teoexact}Let $X$ be an exact $2-$convex rearrangement invariant space.
Then,%
\[
\left\Vert \frac{\sin\pi t}{\pi t}\right\Vert _{X}\leq\left\Vert e^{-\pi
t^{2}/2}\right\Vert _{X}.
\]

\end{theorem}

Let $M$ be an Orlicz function on $[0,\infty)$, i.e., a continuous convex
increasing function such that $M(0)=0$ and $\lim_{u\rightarrow\infty
}M(u)=\infty$. Denote by $L_{M}$ the Orlicz space equipped with the Luxemburg
norm
\begin{equation}
\Vert f\Vert_{L_{M}}:=\inf\Big\{\lambda>0:\,\int_{0}^{\infty}M\left(
\frac{|f(s)|}{\lambda}\right)  \,ds\leq1\Big\}. \label{luxemburg}%
\end{equation}

Applying Theorem \ref{teoexact} to the Orlicz spaces, we obtain

\begin{corollary}
\label{teoexact:cor} If $M(u)$ is an Orlicz function such that $M(\sqrt{u})$
is also an Orlicz function, then
\[
\left\Vert \frac{\sin\pi t}{\pi t}\right\Vert _{L_{M}}\leq\left\Vert e^{-\pi
t^{2}/2}\right\Vert _{L_{M}}.
\]

\end{corollary}

Let us emphasize that it is not possible to prove Theorem \ref{teoexact} using
the \textquotedblleft$L^{p}$ trick\textquotedblright\ of \cite{np} as
described above. Our proof exploits instead the connection between
interpolation theory and the classical theory of majorization (cf. \cite{ca},
\cite{olkin}). We shall now discuss the main features of our approach.

Recall that the classical majorization order $\succeq$ of
Hardy-Littlewood-Polya can be defined in terms of non-increasing
rearrangements as follows\footnote{We refer to the treatise \cite{olkin}
devoted exclusively to majorization and its applications}:%
\begin{equation}
f\succeq g\text{ }\Leftrightarrow\text{ \ }\int_{0}^{t}f^{\ast}(s)\,ds\geq
\int_{0}^{t}g^{\ast}(s)\,ds\text{ for all }t>0. \label{intro2}%
\end{equation}
Another condition, that is often added to the majorization assumptions, is
that $\left\Vert f\right\Vert _{L^{1}}=\left\Vert g\right\Vert _{L^{1}},$
although here often we shall only require an inequality, namely%
\begin{equation}
\left\Vert f\right\Vert _{L^{1}}\geq\left\Vert g\right\Vert _{L^{1}}.
\label{intro2'}%
\end{equation}
One connection between real interpolation and majorization is provided by the
formula for the $K-$functional\footnote{See Section \ref{K-method}} for the
pair $(L^{1},L^{\infty})$ which is given by
\begin{align}
K(t,f;L^{1},L^{\infty})  &  :=\inf\{\left\Vert f_{1}\right\Vert _{L^{1}%
}+t\left\Vert f_{2}\right\Vert _{L^{\infty}}:f=f_{1}+f_{2},f_{1}\in
L^{1},f_{2}\in L^{\infty}\}\nonumber\\
&  =\int_{0}^{t}f^{\ast}(s)\,ds,\;\;t>0 \label{K-f for}%
\end{align}
(cf. \cite[Theorem~5.2.1]{BL} or \cite[Theorem~5.1.6]{BS}).
Thus, (\ref{intro2}) can be written as%
\[
f\succeq g\text{ }\Leftrightarrow\text{ \ }K(t,f;L^{1},L^{\infty})\geq
K(t,g;L^{1},L^{\infty})\text{ for all }t>0.
\]

It is easy to see that (\ref{unou}) (with $C=1$) implies (\ref{intro2}) but
the converse does not hold in general. On the other hand, if $f$ and $g$
satisfy $f\succeq g$, \ then by the Calder\'{o}n-Mityagin theorem (cf.
\cite{ca}, \cite{Mit}, \cite[Theorem 4.6]{BS}), for all rearrangement
invariant function spaces $X$ we have
\[
\left\Vert f\right\Vert _{X}\geq\left\Vert g\right\Vert _{X}.
\]
More generally one can consider other types of majorization associated with
other $K-$, $E-$functionals or some functionals equivalent to them and this is
done in Sections \ref{Int meth} --- \ref{s8} below.

We now present the key argument that leads to our proof of Theorem
\ref{teoexact} (see Section \ref{s7}).

\begin{lemma}
\label{lemma:informal}Suppose that functions $f$ and $g$ satisfy condition
(\ref{distrib}) and, for some $p_{0}\geq1,$ $f\in L^{p_{0}}$ and $\left\Vert
f\right\Vert _{p_{0}}\geq$ $\left\Vert g\right\Vert _{p_{0}}$.
Then, the function $\left\vert g\right\vert ^{p_{0}}$ is majorized by the
function $\left\vert f\right\vert ^{p_{0}}$ in the Hardy-Littlewood-Polya
order,
i.e., for all $t>0,$%
\begin{equation}
\int_{0}^{t}g^{\ast}(s)^{p_{0}}\,ds\leq\int_{0}^{t}f^{\ast}(s)^{p_{0}}\,ds.
\label{achieve}%
\end{equation}

\end{lemma}


\begin{proof}
Since $\lambda_{f}$ and $f^{\ast}$ are generalized inverses for each other,
condition (\ref{distrib}) can be rewritten in terms of rearrangements as
follows: there exists $0<t_{0}<\infty$ such that
\begin{equation}
{f}^{\ast}(t)\left\{
\begin{array}
[c]{cc}%
\geq & {g}^{\ast}(t)\text{ for }t\leq t_{0}\\
\leq & {g}^{\ast}(t)\text{ for }t\geq t_{0}.
\end{array}
\right.  \label{distrib1}%
\end{equation}
To prove (\ref{achieve}) we need to estimate $\int_{0}^{t}g^{\ast}(s)^{p_{0}%
}\,ds$ for all $t>0.$ When $t\leq t_{0},$ (\ref{distrib1}) implies directly
that (\ref{achieve}) holds. It remains to consider the case when $\ t>t_{0}.$
By the cancellation condition $\left\Vert f\right\Vert _{p_{0}}^{p_{0}}\ge$
$\left\Vert g\right\Vert _{p_{0}}^{p_{0}},$ and hence we can write
\begin{align*}
\int_{0}^{t}g^{\ast}(s)^{p_{0}}\,ds  &  =\int_{0}^{\infty}g^{\ast}(s)^{p_{0}%
}\,ds-\int_{t}^{\infty}g^{\ast}(s)^{p_{0}}\,ds\\
&  \le\int_{0}^{\infty}f^{\ast}(s)^{p_{0}}\,ds-\int_{t}^{\infty}g^{\ast
}(s)^{p_{ 0}}\,ds. \label{thisone}%
\end{align*}
Now, for $t>t_{0},$ (\ref{distrib1}) implies that $\int_{t}^{\infty}g^{\ast
}(s)^{p_{0}}\,ds\geq\int_{t}^{\infty}f^{\ast}(s)^{p_{0}}\,ds,$ and, therefore,
we have%
\begin{align*}
\int_{0}^{t}g^{\ast}(s)^{p_{0}}\,ds  &  \leq\int_{0}^{\infty}f^{\ast
}(s)^{p_{0}}\,ds-\int_{t}^{\infty}f^{\ast}(s)^{p_{0}}\,ds\\
&  \leq\int_{0}^{t}f^{\ast}(s)^{p_{0}}\,ds,\text{ for all }t>0,
\end{align*}
as we wished to show.
\end{proof}

In Section \ref{differential inequalities}, which is somewhat more informal,
we consider some comparison of norms that is connected with developments
associated with derivatives of norms and entropies (e.g. Log Sobolev
inequalities, Commutator Inequalities, etc.).

Finally, we have added an Appendix, addressed primarily to those readers who
might not be familiar with interpolation theory. In Section \ref{s9}, we give
a simplified version of the Nazarov-Podkorytov Lemma, showing that an
inequality of the type $\Vert g\Vert_{p}\leq\Vert f\Vert_{p},$ for $p$ large
enough, can be obtained whenever we know that $g^{\ast}(t)\leq f^{\ast}(t),$
for all small $t>0$ (see Proposition \ref{cor1d}). However, in contrast to
\cite{np}, in this regime we cannot indicate an initial value $p_{0},$ such
the norm inequality holds for $p\geq p_{0}.$ We close the paper (cf. Section
\ref{s10}) showing a close connection between the Nazarov-Podkorytov Lemma and
the well-known classical Karamata inequality. In particular, general
Nazarov-Podkorytov type results, for modular forms, follow.

\section{Results of Nazarov-Podkorytov type for Lorentz $L(p,q)$%
-spaces\label{s2}}

\subsection{Sufficient conditions\label{s2a}}

The main purpose of this section is to formulate a sharp result of
Nazarov-Podkorytov type in the context of the $L(p,q)$-spaces. Although we
retain the basic idea of the original argument, the extension is not
completely straightforward and, in particular, requires suitable restrictions
on the parameters. In Section \ref{s3} we prove the necessity of the
aforementioned conditions.

\begin{proposition}
\label{cor1b} Let $0<p_{0},q_{0},p,q<\infty$. Suppose that $f$ and $g$ are
measurable functions such that (\ref{distrib}) holds for some $\tau_{0}>0$.
Moreover, suppose that one of the following conditions is satisfied: (a)
$q>q_{0}$, $q/p=q_{0}/p_{0}$, or (b) $q=q_{0}$, $p>p_{0}$.

Furthermore, suppose that%
\begin{equation}
\int_{0}^{\infty}(f^{\ast}(t)^{q_{0}}-g^{\ast}(t)^{q_{0}})\,d(t^{q_{0}/p_{0}%
})\geq0 \label{penso(pr)}%
\end{equation}
and, moreover,
\begin{equation}
(f^{\ast}(t)^{q}-g^{\ast}(t)^{q})t^{q/p-1}\in L_{1}(0,\infty).
\label{penso(pra)}%
\end{equation}
Then, we have
\begin{equation}
\int_{0}^{\infty}(f^{\ast}(t)^{q}-g^{\ast}(t)^{q})\,d(t^{q/p})\geq0.
\label{penso2(pr)}%
\end{equation}
As a consequence, if $f,g\in L(p,q),$ then
\[
\Vert g\Vert_{L(p,q)}\leq\Vert f\Vert_{L(p,q)}.
\]

\end{proposition}

\begin{proof}
The argument follows the main lines of the proof of the Nazarov-Podkorytov
Lemma, as explained in the Introduction.

\textbf{(a)} First, we prove that condition \eqref{penso(pra)} ensures that
\begin{equation}
\int_{0}^{\infty}(f^{\ast}(t)^{q}-g^{\ast}(t)^{q})t^{q/p-1}\,dt=p\int%
_{0}^{\infty}\left(  (\lambda_{f}(u))^{q/p}-(\lambda_{g}(u))^{q/p}\right)
u^{q-1}\,du. \label{penso2(prb)}%
\end{equation}

We set $h_{1}(t):=f^{\ast}(t)^{q}$, $h_{2}(t):=g^{\ast}(t)^{q}$, and
$h(t):=\min_{i=1,2}h_{i}(t)$. Then, $h_{1}-h\geq0$ and $h_{1}-h\leq
|h_{1}-h_{2}|$, whence from \eqref{penso(pra)} it follows that $(h_{1}%
(t)-h(t))t^{q/p-1}\in L_{1}(0,\infty)$. Let
\[
A_{i}:=\{(t,v)\in(0,\infty)\times(0,\infty):\,h(t)\leq v\leq h_{i}%
(t)\},\;i=1,2.
\]
Since $h_{1}$ and $h_{2}$ are decreasing, we get
\[
A_{i}=\{(t,v)\in(0,\infty)\times(0,\infty):\,\lambda_{h}(v)\leq t\leq
\lambda_{h_{i}}(v)\},\;i=1,2.
\]
Computing iterated yields for $i=1,2$
\[
\int_{0}^{\infty}\int_{0}^{\infty}\chi_{A_{i}}(t,v)t^{q/p-1}\,dv\,dt=\int%
_{0}^{\infty}(h_{i}(t)-h(t))t^{q/p-1}\,dt<\infty
\]
and
\begin{align*}
\int_{0}^{\infty}\int_{0}^{\infty}\chi_{A_{i}}(t,v)t^{q/p-1}\,dt\,dv  &
=\int_{0}^{\infty}\int_{\lambda_{h}(v)}^{\lambda_{h_{i}}(v)}t^{q/p-1}%
\,dt\,dv\\
&  =\frac{p}{q}\int_{0}^{\infty}\left(  \lambda_{h_{i}}(v)^{q/p}-\lambda
_{h}(v)^{q/p}\right)  \,dv.
\end{align*}
Therefore, applying Fubini's theorem, we obtain for $i=1,2$
\[
\int_{0}^{\infty}(h_{i}(t)-h(t))t^{q/p-1}\,dt=\frac{p}{q}\int_{0}^{\infty
}\left(  \lambda_{h_{i}}(v)^{q/p}-\lambda_{h}(v)^{q/p}\right)  \,dv,
\]
whence
\[
\int_{0}^{\infty}(h_{2}(t)-h_{1}(t))t^{q/p-1}\,dt=\frac{p}{q}\int_{0}^{\infty
}\left(  \lambda_{h_{2}}(v)^{q/p}-\lambda_{h_{1}}(v)^{q/p}\right)  \,dv.
\]
From the definitions, $\lambda_{h_{1}}(v)=\lambda_{f}(v^{1/q})$,
$\lambda_{h_{2}}(v)=\lambda_{g}(v^{1/q})$; therefore, after a suitable change
of variables, we arrive at \eqref{penso2(prb)}. Introducing the parameter
$\tau_{0}$, we obtain
\begin{equation}
\int_{0}^{\infty}(f^{\ast}(t)^{q}-g^{\ast}(t)^{q})\,d(t^{q/p})=q\tau_{0}%
^{q-1}\int_{0}^{\infty}\left(  \lambda_{f}(u)^{q/p}-\lambda_{g}(u)^{q/p}%
\right)  \left(  \frac{u}{\tau_{0}}\right)  ^{q-1}du. \label{penso3(prbb)}%
\end{equation}

For $0<\gamma,\beta<\infty,$ let
\[
\psi_{f,g}(\gamma,\beta):=\int_{0}^{\infty}\left(  \lambda_{f}(u)^{\beta
/\gamma}-\lambda_{g}(u)^{\beta/\gamma}\right)  \left(  \frac{u}{\tau_{0}%
}\right)  ^{\beta-1}\,du.
\]
Using the fact that $q/p=q_{0}/p_{0}$, we find
\[
\psi_{f,g}(p,q)-\psi_{f,g}(p_{0},q_{0})=\int_{0}^{\infty}\left(  \left(
\frac{u}{\tau_{0}}\right)  ^{q-1}-\left(  \frac{u}{\tau_{0}}\right)
^{q_{0}-1}\right)  (\lambda_{f}(u)^{q_{0}/p_{0}}-\lambda_{g}(u)^{q_{0}/p_{0}%
})\,du.
\]
Since $q>q_{0}$, from (\ref{distrib}) it is easy to verify that the factors
\[
\left(  \lambda_{f}(u)^{q_{0}/p_{0}}-\lambda_{g}(u)^{q_{0}/p_{0}}\right)
\;\;\mbox{and}\;\;\left(  \left(  \frac{u}{\tau_{0}}\right)  ^{q-1}-\left(
\frac{u}{\tau_{0}}\right)  ^{q_{0}-1}\right)
\]
have the same signs on each of the intervals $(0,\tau_{0})$ and $(\tau
_{0},\infty).$ Thus,
\[
\psi_{f,g}(p,q)-\psi_{f,g}(p_{0},q_{0})\geq0.
\]
Consequently, by (\ref{penso(pr)}), \eqref{penso3(prbb)} and the fact that
$\psi_{f,g}(p_{0},q_{0})\geq0$, we see that $\psi_{f,g}(p,q)\geq0$, and
(\ref{penso2(pr)}) follows.

\textbf{(b)} Since $h^{\ast}$ and $\lambda_{h}$ are generalized inverses of
each other, from (\ref{distrib}) it follows that there exists $t_{0}$ such
that \eqref{distrib1} holds.
Then, we have
\[
\int_{0}^{\infty}(f^{\ast}(t)^{q}-g^{\ast}(t)^{q})\,d(t^{q/p})=\frac{q}%
{p}t_{0}^{q/p-1}\int_{0}^{\infty}\left(  f^{\ast}(t)^{q}-g^{\ast}%
(t)^{q}\right)  \left(  \frac{t}{t_{0}}\right)  ^{q/p-1}\,dt.
\]
Next, setting
\[
\kappa_{f,g}(\gamma,\beta):=\int_{0}^{\infty}\left(  f^{\ast}(t)^{\beta
}-g^{\ast}(t)^{\beta}\right)  \left(  \frac{t}{t_{0}}\right)  ^{\beta
/\gamma-1}\,dt,
\]
we verify by computation that
\[
\kappa_{f,g}(p,q_{0})-\kappa_{f,g}(p_{0},q_{0})=\int_{0}^{\infty}\left(
f^{\ast}(t)^{q_{0}}-g^{\ast}(t)^{q_{0}}\right)  \left(  \left(  \frac{t}%
{t_{0}}\right)  ^{q_{0}/p-1}-\left(  \frac{t}{t_{0}}\right)  ^{q_{0}/p_{0}%
-1}\right)  \,dt.
\]
Using (\ref{distrib1}) and the fact that $p>p_{0}$, we can analyze the sign of
the integrand, as in the proof of case \textbf{(a) }above, to conclude that
\[
\kappa_{f,g}(p,q_{0})-\kappa_{f,g}(p_{0},q_{0})\geq0.
\]
Appealing once again to (\ref{penso(pr)}), it follows that $\kappa
_{f,g}(p,q_{0})\geq0$, and (\ref{penso2(pr)}) follows.
\end{proof}

\medskip\vskip0.3cm

The next result shows that if the supports of $f$ and $g$ have finite
measure\footnote{By abuse of language we will say simply that $f$ and $g$ have
\textquotedblleft finite support\textquotedblright.}, then condition
\eqref{distrib} can be somewhat weakened.

\begin{proposition}
\label{prop1} Let $0<p_{0},p<\infty$, $1<q_{0}<q<\infty,$ $q/p=q_{0}/p_{0}$,
and let $f$, $g$ be measurable functions with finite support satisfying
\eqref{penso(pr)}. Moreover, suppose that $f,g\in L(p,q)$ and there exists
$\tau_{1}\in\lbrack0,\infty)$ such that
\begin{equation}
\int_{\tau}^{\infty}\lambda_{g}(s)^{q_{0}/p_{0}}\,ds\leq\int_{\tau}^{\infty
}\lambda_{f}(s)^{q_{0}/p_{0}}\,ds\text{ if }\tau\geq\tau_{1}, \label{eq104a}%
\end{equation}
and
\begin{equation}
\int_{\tau}^{\infty}\lambda_{f}(s)^{q_{0}/p_{0}}\,ds\leq\int_{\tau}^{\infty
}\lambda_{g}(s)^{q_{0}/p_{0}}\,ds\text{ if }\tau\leq\tau_{1}. \label{eq105a}%
\end{equation}
Then,
\[
\Vert g\Vert_{L(p,q)}\leq\Vert f\Vert_{L(p,q)}.
\]

\end{proposition}

\begin{proof}
We only consider in detail the case $\tau_{1}>0$, the case when $\tau_{1}=0$
is easier.

Using \eqref{quasi-norm} and integration by parts, we can write
\begin{align*}
\int_{0}^{\infty}(f^{\ast}(t)^{q}-g^{\ast}(t)^{q})\,d(t^{q/p})  &  =q\tau
_{1}^{q-1}\int_{0}^{\infty}\left(  \lambda_{f}(s)^{q/p}-\lambda_{g}%
(s)^{q/p}\right)  \left(  \frac{s}{\tau_{1}}\right)  ^{q-1}ds\\
&  =-q\tau_{1}^{q-1}\int_{0}^{\infty}\left(  \frac{s}{\tau_{1}}\right)
^{q-1}\,d\left(  \int_{s}^{\infty}\left(  \lambda_{f}(u)^{q/p}-\lambda
_{g}(u)^{q/p}\right)  \,du\right) \\
&  =q(q-1)\tau_{1}^{q-2}\int_{0}^{\infty}\left(  \frac{s}{\tau_{1}}\right)
^{q-2}\int_{s}^{\infty}\left(  \lambda_{f}(u))^{q/p}-(\lambda_{g}%
(u)^{q/p}\right)  \,du\,ds\\
&  +(I)+(II),
\end{align*}
where
\[
(I)=q\lim_{s\rightarrow0}s^{q-1}\int_{s}^{\infty}\left(  \lambda_{f}%
(u)^{q/p}-\lambda_{g}(u)^{q/p}\right)  \,du,
\]%
\[
(II)=-q\lim_{s\rightarrow\infty}s^{q-1}\int_{s}^{\infty}\left(  \lambda
_{f}(u)^{q/p}-\lambda_{g}(u)^{q/p}\right)  \,du.
\]
We shall now verify that $(I)=0$ and $(II)=0.$ The first equation will be
proved once we show that
\[
(A):=\lim_{s\rightarrow0}s^{q-1}\int_{s}^{\infty}\lambda_{f}(u)^{q/p}%
\,du=0,\;\;(B):=\lim_{s\rightarrow0}s^{q-1}\int_{s}^{\infty}\lambda
_{g}(u)^{q/p}\,du=0.
\]
Let us note that, as it will become apparent, the proof that $(A)=0$ can be
applied verbatim to verify that $(B)=0$. Let $\left\Vert f\right\Vert _{0}$
denote the measure of the support of $f$ (see \eqref{e1a}), then since $f\in
L(p,q)$ and $q>1$, we obtain
\begin{align*}
(A)  &  =\lim_{s\rightarrow0}s^{q-1}\left\{  \int_{s}^{1}\lambda_{f}%
(u)^{q/p}\,du+\int_{1}^{\infty}\lambda_{f}(u)^{q/p}\,du\right\} \\
&  \leq\lim_{s\rightarrow0}s^{q-1}\left\{  \left\Vert f\right\Vert _{0}%
^{q/p}+\int_{1}^{\infty}\lambda_{f}(u)^{q/p}u^{q-1}\,du\right\} \\
&  \leq\lim_{s\rightarrow0}s^{q-1}\left(  \left\Vert f\right\Vert _{0}%
^{q/p}+\left\Vert f\right\Vert _{L(p,q)}^{q}\right)  =0.
\end{align*}
The fact that $(II)=0$ follows by observing that, for $s>0,$%
\[
s^{q-1}\int_{s}^{\infty}\lambda_{f}(u)^{q/p}du\leq\int_{s}^{\infty}\lambda
_{f}(u)^{q/p}u^{q-1}\,du,
\]
and that, in view of the fact that $f\in L(p,q),$ the right-hand side tends to
zero when $s\rightarrow\infty.$

Summarizing the previous discussion, we have shown that integrating by parts
yields
\[
\int_{0}^{\infty}(f^{\ast}(t)^{q}-g^{\ast}(t)^{q})\,d(t^{q/p})=q(q-1)\tau
_{1}^{q-2}\int_{0}^{\infty}\left(  \frac{s}{\tau_{1}}\right)  ^{q-2}\int%
_{s}^{\infty}\left(  \lambda_{f}(u)^{q/p}-\lambda_{g}(u)^{q/p}\right)
\,du\,ds.
\]
Let
\[
\Psi_{f,g}(p,q):=\int_{0}^{\infty}\left(  \frac{s}{\tau_{1}}\right)
^{q-2}\int_{s}^{\infty}\left(  \lambda_{f}(u)^{q/p}-\lambda_{g}(u)^{q/p}%
\right)  \,du\,ds.
\]
On account of the fact that $q/p=q_{0}/p_{0}$, we can write
\[
\Psi_{f,g}(p,q)-\Psi_{f,g}(p_{0},q_{0})=\int_{0}^{\infty}\left(  \left(
\frac{s}{\tau_{1}}\right)  ^{q-2}-\left(  \frac{s}{\tau_{1}}\right)
^{q_{0}-2}\right)  \int_{s}^{\infty}\left(  \lambda_{f}(u)^{q_{0}/p_{0}%
}-\lambda_{g}(u)^{q_{0}/p_{0}}\right)  \,du\,ds.
\]
Using inequalities \eqref{eq104a} and \eqref{eq105a}, we see that the factors
\[
\int_{s}^{\infty}\left(  \lambda_{f}(u)^{q_{0}/p_{0}}-\lambda_{g}%
(u)^{q_{0}/p_{0}}\right)  \,du\;\;\mbox{and}\;\;\left(  \left(  \frac{s}%
{\tau_{1}}\right)  ^{q-2}-\left(  \frac{s}{\tau_{1}}\right)  ^{q_{0}%
-2}\right)
\]
have the same signs on each of the intervals $(0,\tau_{1})$ and $(\tau
_{1},\infty)$. Therefore,
\[
\Psi_{f,g}(p,q)-\Psi_{f,g}(p_{0},q_{0})\geq0,
\]
and, since $\Psi_{f,g}(p_{0},q_{0})\geq0$ by \eqref{penso(pr)}, we see that
$\Psi_{f,g}(p,q)\geq0$, which is equivalent to \eqref{penso2(pr)}.
\end{proof}

Applying Proposition \ref{prop1} with $p_{0}=q_{0}>1$, we derive the following
version of the Nazarov-Podkorytov Lemma in the $L^{p}$-setting.

\begin{corollary}
\label{cor1a} Let $p\ge p_{0}>1$. Suppose that functions $f\in L^{p}$, $g\in
L^{p}$ have finite support and such that $\Vert g\Vert_{p_{0}}\leq\Vert
f\Vert_{p_{0}}$. Suppose, moreover, that there exists $\tau_{1}\in
\lbrack0,\infty)$ such that
\begin{equation}
\int_{\tau}^{\infty}\lambda_{g}(s)\,ds\leq\int_{\tau}^{\infty}\lambda
_{f}(s)\,ds\text{ if }\tau\geq\tau_{1}, \label{eq104b}%
\end{equation}%
\begin{equation}
\int_{\tau}^{\infty}\lambda_{f}(s)\,ds\leq\int_{\tau}^{\infty}\lambda
_{g}(s)\,ds\text{ if }\tau\leq\tau_{1}. \label{eq105b}%
\end{equation}
Then, $\Vert g\Vert_{p}\leq\Vert f\Vert_{p}$.
\end{corollary}

\begin{remark}
One can easily see that conditions (\ref{eq104b}) and (\ref{eq105b}) are
weaker than condition (\ref{distrib}), that is, (\ref{distrib}) implies
(\ref{eq104b}) and (\ref{eq105b}) for some $\tau_{1}\in[0,\tau_{0}]$, but the
converse does not hold in general.


\end{remark}

\subsection{Sharp conditions\label{s3}}

Next we show that Proposition \ref{cor1b} implies a more general result of the
Nazarov-Podkorytov type for $L(p,q)$-spaces.

\begin{corollary}
\label{cor1ab} Let $0<p_{0}\leq p<\infty$, $0<q_{0}\leq q<\infty$ and $q/p\leq
q_{0}/p_{0}$. Suppose that $f$ and $g$ are measurable functions such that
(\ref{distrib}) holds for some $\tau_{0}>0$. Moreover, assume that
\[
\int_{0}^{\infty}(f^{\ast}(t)^{q_{0}}-g^{\ast}(t)^{q_{0}})\,d(t^{q_{0}/p_{0}%
})\geq0
\]
and
\[
(f^{\ast}(t)^{q}-g^{\ast}(t)^{q})t^{q/p-1}\in L^{1}(0,\infty).
\]
Then, we have
\[
\int_{0}^{\infty}(f^{\ast}(t)^{q}-g^{\ast}(t)^{q})\,d(t^{q/p})\geq0.
\]
Consequently, if $f,g\in L(p,q),$ then
\[
\Vert g\Vert_{L(p,q)}\leq\Vert f\Vert_{L(p,q)}.
\]

\end{corollary}

Before going through the proof of Corollary \ref{cor1ab}, it will be
convenient to adopt the following notation.

Let $0<{p_{0}},{q_{0}}<\infty$ be fixed. Denote by ${\mathcal{A}}({p_{0}%
},{q_{0}})$ the set of all pairs of indices $({p},{q})$, $0<{p},{q}<\infty,$
such that for arbitrary measurable functions $f,g$, that satisfy condition
\eqref{distrib} for some $\tau_{0}>0$ and the inequality
\[
\int_{0}^{\infty}(f^{\ast}(t)^{{q_{0}}}-g^{\ast}(t)^{{q_{0}}})\,d(t^{{q_{0}%
}/{p_{0}}})\geq0,
\]
from the validity of the condition
\[
(f^{\ast}(t)^{{q}}-g^{\ast}(t)^{{q}})t^{{q}/{p}-1}\in L^{1}(0,\infty)
\]
we can deduce that
\[
\int_{0}^{\infty}(f^{\ast}(t)^{{q}}-g^{\ast}(t)^{{q}})\,d(t^{{q}/{p}})\geq0.
\]

The main aim of this section is to provide a complete characterization of the
set ${\mathcal{A}}({p_{0}},{q_{0}})$ showing the sharpness of the conditions
imposed on parameters of Lorentz spaces $L(p,q)$ in Corollary \ref{cor1ab}.

For our purposes the following transitivity property, which is a routine
unraveling of the definitions, will be useful.

\begin{lemma}
\label{lemma trans1} If $(p_{1},q_{1})\in{\mathcal{A}}(p,q)$ and $(p_{2}%
,q_{2})\in{\mathcal{A}}(p_{1},q_{1})$, then $(p_{2},q_{2})\in{\mathcal{A}%
}(p,q)$.
\end{lemma}

\begin{proof}
[Proof of Corollary \ref{cor1ab}]Since ${q}/{p}\leq q_{0}/p_{0}$, we can
select $\gamma\geq p_{0}$ such that ${q}=q_{0}{p}/\gamma$. Then, from the
definition of the set ${\mathcal{A}}(p,q)$ and Proposition \ref{cor1b} it
follows that $(\gamma,q_{0})\in{\mathcal{A}}(p_{0},q_{0})$ and $({p},{q}%
)\in{\mathcal{A}}(\gamma,q_{0})$. Therefore, by Lemma \ref{lemma trans1}, we
have $({p},{q})\in{\mathcal{A}}(p_{0},q_{)})$. This completes the proof.
\end{proof}


\begin{theorem}
\label{prop of Dlta} Let $0<{p_{0}},{q_{0}}<\infty,$ then%
\begin{equation}
{\mathcal{A}}({p_{0}},{q_{0}})=\{({p},{q}):\,{p}\geq{p_{0}},{q}\geq{q_{0}}%
,{q}/{p}\leq{q_{0}}/{p_{0}}\}. \label{righthandside}%
\end{equation}

\end{theorem}

\begin{proof}
For notational convenience, let $\Delta({p_{0}},{q_{0}})$ be the set of pairs
$({p},{q})$ defined by the right-hand side of (\ref{righthandside}).

By Corollary \ref{cor1ab}, we have $\Delta({p_{0}},{q_{0}})\subset
{\mathcal{A}}({p_{0}},{q_{0}})$.
To prove the converse, it suffices to show that from $({p},{q})\not \in
\Delta({p_{0}},{q_{0}})$ it follows $({p},{q})\not \in {\mathcal{A}}({p_{0}%
},{q_{0}})$. We consider three cases separately.

(i) ${p}<{p_{0}}$. Let
\begin{equation}
f(t)=\chi_{\lbrack0,1]}(t),\;g(t)=a\chi_{\lbrack0,1]}(t)+b\chi_{(1,h]}%
(t),\;\text{ where }\;1>a\geq b>0,h>1. \label{asbefore}%
\end{equation}
By inspection, we see that $f(t)>g(t)$ if $0<t\leq1$, and $f(t)\leq g(t)$ if
$t>1.$ Moreover, it is also plain that $f(t)=f^{\ast}(t)$, $g(t)=g^{\ast}(t)$.
A simple computation yields
\[
\Vert f\Vert_{r,s}^{s}=1\;\;\mbox{and}\;\;\Vert g\Vert_{r,s}^{s}=a^{s}%
+b^{s}(h^{s/r}-1)\;\;\mbox{for all}\;\;0<r,s<\infty.
\]
Hence, from
\begin{equation}
a^{{q_{0}}}+b^{{q_{0}}}(h^{{q_{0}}/{p_{0}}}-1)=1, \label{norm}%
\end{equation}
it follows that $\Vert f\Vert_{{p_{0}},{q_{0}}}=\Vert g\Vert_{{p_{0}},{q_{0}}%
}$. We now show that for a suitable choice of the parameters $a,b,h$ which
satisfy \eqref{norm} we have $\Vert g\Vert_{{p},{q}}>\Vert f\Vert_{{p},{q}}$.
Since $\Vert f\Vert_{{r},{s}}=1$ for all $0<r,s<\infty$, it suffices to prove
that, under condition \eqref{norm}, it holds
\begin{equation}
a^{{q}}+b^{{q}}(h^{{q}/{p}}-1)>1. \label{aim}%
\end{equation}

Select $a=b$. Then solving \eqref{norm} for $a$, we get $a=b=h^{-1/{p_{0}}}$.
Thus, computing the left-hand side of (\ref{aim}), and taking into account the
fact that ${p}<{p_{0}}$ and $h>1,$ yields
\[
h^{-{q}/{p_{0}}}+h^{-{q}/{p_{0}}}(h^{{q}/{p}}-1)=h^{{q}(1/{p}-1/{p_{0}})}>1,
\]
and (\ref{aim}) follows.

(ii) ${q}<{q_{0}}$. Let $f$ and $g$ be defined once again by (\ref{asbefore}).
Since for any $\alpha>0$
\[
\lim_{h\rightarrow1}\frac{h^{\alpha}-1}{h-1}=\alpha,
\]
we can choose $h>1$ so that
\begin{equation}
h^{{q}/{p}}-1>\frac{{q}}{2{p}}(h-1)\;\;\mbox{and}\;\;h^{{q_{0}}/{p_{0}}%
}-1<\frac{2{q_{0}}}{{p_{0}}}(h-1). \label{normB}%
\end{equation}
Fix $h>1$ satisfying \eqref{normB}. Using the fact that ${q}<{q_{0}}$, select
$b\in(0,h^{-1/{p_{0}}})$, such that
\begin{equation}
b^{{q}-{q_{0}}}>\frac{4{q_{0}}{p}}{{p_{0}}{q}}. \label{normA}%
\end{equation}
Setting $a=(1-b^{{q_{0}}}(h^{{q_{0}}/{p_{0}}}-1))^{1/{q_{0}}}$, we see that
$0<b<a<1,$ and the triple $(a,b,h)$ satisfies \eqref{norm}.

Now we show that (\ref{aim}) holds. Indeed, taking into account \eqref{norm},
(\ref{normA}), (\ref{normB}) and the fact that ${q}<{q_{0}}$, we can write
\begin{align*}
a^{{q}}+b^{{q}}(h^{{q}/{p}}-1)-1  &  =a^{{q}}-a^{{q_{0}}}+b^{{q}}(h^{{q}/{p}%
}-1)-1+a^{{q_{0}}}\\
&  =(a^{{q}}-a^{{q_{0}}})+(b^{{q}}(h^{{q}/{p}}-1)-b^{{q_{0}}}(h^{{q_{0}%
}/{p_{0}}}-1))\\
&  >\Big(\frac{{q}}{2{p}}b^{{q}}-\frac{2{q_{0}}}{{p_{0}}}b^{{q_{0}}%
}\Big)(h-1)\\
&  =b^{{q_{0}}}\Big(\frac{{q}}{2{p}}b^{{q}-{q_{0}}}-\frac{2{q_{0}}}{{p_{0}}%
}\Big)(h-1)\\
&  >0.
\end{align*}

(iii) ${q}/{p}>{q_{0}}/{p_{0}}$. Let $a,b,h$ be parameters such that
$a>1>b>0$, $0<h<1,$ whose precise values will be specified later. Let
\[
f(t)=a\chi_{\lbrack0,h]}(t)+b\chi_{(h,1]}(t),\;\;g(t)=\chi_{\lbrack0,1]}(t),
\]
It is plain that $f(t)>g(t)$ if $0<t\leq h$, while $f(t)<g(t)$ if $t>h$.
Moreover, since $f,$ and $g$, are decreasing, $f(t)=f^{\ast}(t)$,
$g(t)=g^{\ast}(t)$. Our aim is to prove that for a suitable choice of the
parameters $a,b,h$ we have that
\begin{equation}
\Vert f\Vert_{{p_{0}},{q_{0}}}=\Vert g\Vert_{{p_{0}},{q_{0}}}%
\;\;\mbox{and}\;\;\Vert f\Vert_{{p},{q}}<\Vert g\Vert_{{p},{q}}. \label{normC}%
\end{equation}
By computation,
\[
\Vert f\Vert_{{p_{0}},{q_{0}}}^{{q_{0}}}=a^{{q_{0}}}h^{{q_{0}}/{p_{0}}%
}+b^{{q_{0}}}(1-h^{{q_{0}}/{p_{0}}})=b^{{q_{0}}}+(a^{{q_{0}}}-b^{{q_{0}}%
})h^{{q_{0}}/{p_{0}}}\;\;\mbox{and}\;\;\Vert g\Vert_{{p_{0}},{q_{0}}}^{{q_{0}%
}}=1.
\]
Therefore, $\Vert f\Vert_{{p_{0}},{q_{0}}}=\Vert g\Vert_{{p_{0}},{q_{0}}}$ is
equivalent to
\begin{equation}
b^{{q_{0}}}+(a^{{q_{0}}}-b^{{q_{0}}})h^{{q_{0}}/{p_{0}}}=1. \label{norm1}%
\end{equation}
Likewise, $\Vert f\Vert_{{p},{q}}<\Vert g\Vert_{{p},{q}}$ can be rewritten as
\[
a^{{q}}h^{{q}/{p}}+b^{{q}}(1-h^{{q}/{p}})<1=b^{{q_{0}}}+(a^{{q_{0}}}%
-b^{{q_{0}}})h^{{q_{0}}/{p_{0}}},
\]
and after some further elementary manipulations, can be seen to be equivalent
to
\begin{equation}
(a^{{q_{0}}}-b^{{q_{0}}})h^{{q_{0}}/{p_{0}}-{q}/{p}}-(a^{{q}}-b^{{q}}%
)>(b^{{q}}-b^{{q_{0}}})h^{-{q}/{p}}. \label{norm2}%
\end{equation}

Clearly, we can additionally assume here that ${q_{0}}\leq{q}$. Therefore,
from the inequality $b<1$ it follows that the right-hand side of (\ref{norm2})
is non-positive. Consequently, (\ref{norm2}) would be proved once we show
that
\begin{equation}
(a^{{q_{0}}}-b^{{q_{0}}})h^{{q_{0}}/{p_{0}}-{q}/{p}}-(a^{{q}}-b^{{q}})>0.
\label{normD}%
\end{equation}
Since from (\ref{norm1}) it follows that
\begin{equation}
h=\left(  \frac{1-b^{{q_{0}}}}{a^{{q_{0}}}-b^{{q_{0}}}}\right)  ^{{p_{0}%
}/{q_{0}}}, \label{normE}%
\end{equation}
if we let $\beta={{q}{p_{0}}}/({p}{q_{0}})$, we see that \eqref{normD} is
equivalent to the inequality
\[
\left(  \frac{1-b^{{q_{0}}}}{a^{{q_{0}}}-b^{{q_{0}}}}\right)  ^{1-\beta}%
>\frac{a^{{q}}-b^{{q}}}{a^{{q_{0}}}-b^{{q_{0}}}},
\]
or
\begin{equation}
(1-b^{{q_{0}}})^{1-\beta}>(a^{{q}}-b^{{q}})(a^{{q_{0}}}-b^{{q_{0}}})^{-\beta}.
\label{ultima}%
\end{equation}
Let $a>1$ be fixed. Then, since $\beta>1$, we can choose $b<1$ sufficiently
close to $1$ so that (\ref{ultima}) holds. Let $h\in(0,1)$ be defined by
\eqref{normE}. Then the triple $(a,b,h)$ satisfies both \eqref{norm1} and
\eqref{norm2}. Consequently, \eqref{normC} holds completing the proof of the theorem.
\end{proof}

\subsection{Integral inequalities of Ball's type for $L(p,q)$-norms}

\label{sec:integralineq}

As we have seen in the Introduction, Ball's integral inequality \eqref{intro1}
can be reformulated as the comparison of the $L^{p}(0,\infty)$-norms of the
functions $f$ and $g$ defined by \eqref{f&g}. Since for $p=2$ inequality
\eqref{tres} becomes an identity and the distribution functions of $f$ and $g$
satisfy \eqref{distrib} for some $\tau_{0}>0$, we can apply Corollary
\ref{cor1ab} to obtain the following $L(p,q)$-version of Ball's integral inequality.

\begin{theorem}
\label{teor<} For every $2\leq{q}\leq{p}<\infty,$ we have
\begin{equation}
\int_{0}^{+\infty}\left(  \left(  \frac{\sin\pi t}{\pi t}\right)  ^{\ast
}\right)  ^{{q}}t^{{q}/{p}-1}\,dt\leq\int_{0}^{+\infty}e^{-\pi{q}t^{2}%
/2}t^{{q}/{p}-1}\,dt. \label{Ball1}%
\end{equation}
Equivalently, for all $q\geq2$ and $0<\alpha\leq1,$
\begin{equation}
\int_{0}^{+\infty}\left(  \left(  \frac{\sin t}{t}\right)  ^{\ast}\right)
^{{q}}t^{\alpha-1}\,dt\leq\frac{1}{2}q^{-\alpha/2}(2\pi)^{\alpha/2}%
\Gamma(\alpha/2), \label{Ball2}%
\end{equation}
where $\Gamma(x)$ is the gamma-function.
\end{theorem}

\begin{proof}
Inequality \eqref{Ball1} is an immediate consequence of Corollary \ref{cor1ab}
with $p_{0}=q_{0}=2$.
To derive \eqref{Ball2}, we make a change of variables: $u=\pi t$ on the
left-hand side of \eqref{Ball1}, $u=\pi qt^{2}/2$ on the right-hand side, and
let $\alpha=q/p$. Then, by straightforward calculations, we obtain
\begin{align*}
\int_{0}^{+\infty}\left(  \left(  \frac{\sin u}{u}\right)  ^{\ast}\right)
^{{q}}u^{\alpha-1}\,du  &  \leq\frac{1}{2}q^{-\alpha/2}(2\pi)^{\alpha/2}%
\int_{0}^{+\infty}e^{-u}u^{\alpha/2-1}\,du\\
&  =\frac{1}{2}q^{-\alpha/2}(2\pi)^{\alpha/2}\Gamma(\alpha/2).
\end{align*}

\end{proof}

Recently, in \cite{CKT-20}, Chasapis, K\"{o}nig and Tkocz have obtained a
sharp $L_{-1}-L_{2}$ Khintchine type inequality (see \cite{Kh}, \cite{Haag}),
based on an extension of Ball's integral inequality they proved in the same
paper. It will be instructive to compare the version of Ball's inequality
given in \cite{CKT-20} with our inequality \eqref{Ball2}. For this purpose,
let us define for every $\alpha\in(0,1)$ the constants
Let us define first for every $\alpha\in(0,1)$ the constants
\[
c_{2}(\alpha):=\frac{2^{1-\alpha/2}}{(1-\alpha)(2-\alpha)},\;\;c_{\infty
}(\alpha):=\frac{(3/2)^{\alpha/2}}{\sqrt{\pi}}\Gamma\left(  \frac{1-\alpha}%
{2}\right)  ,\;\;C_{\alpha}:=\max\{c_{2}(\alpha),c_{\infty}(\alpha)\},
\]
which are of a certain probabilistic nature (cf. \cite{CKT-20}). It was shown
in \cite{CKT-20}, that there exists a unique $\alpha_{0}=0.793...$ such that
$c_{2}(\alpha_{0})=c_{\infty}(\alpha_{0})$ with $C_{\alpha}=c_{\infty}%
(\alpha)$ if $\alpha\in(0,\alpha_{0})$ and $C_{\alpha}=c_{2}(\alpha)$ if
$\alpha\in(\alpha_{0},1)$. Then (cf. \cite{CKT-20}), for all $q\geq2$ and
$0<\alpha<1,$
\begin{equation}
\int_{0}^{+\infty}\left\vert \frac{\sin t}{t}\right\vert ^{{q}}t^{\alpha
-1}\,dt\leq2^{\alpha-1}q^{-\alpha/2}\sqrt{\pi}\frac{\Gamma(\alpha/{2})}%
{\Gamma((1-\alpha)/{2})}C_{\alpha}. \label{Ball3}%
\end{equation}

Since the function $t^{\alpha-1}$ decreases for $t>0$, the term on the
left-hand side of inequality \eqref{Ball2} is larger than the corresponding
term of inequality \eqref{Ball3} (see e.g. \cite[\S \,II.2]{KPS}). Let us now
compare the constants from the right-hand sides of these inequalities.

Suppose first $\alpha\in(0,\alpha_{0})$. Then, $C_{\alpha}=c_{\infty}(\alpha)$
and the right-hand side of \eqref{Ball3} is equal to
\[
2^{\alpha-1}q^{-\alpha/2}(3/2)^{\alpha/2}\Gamma(\alpha/{2})=\frac{1}%
{2}q^{-\alpha/2}6^{\alpha/2}\Gamma(\alpha/{2}).
\]
We see that the ratio of the constant from inequality \eqref{Ball2} to that of
inequality \eqref{Ball3} is $(\pi/3)^{\alpha/2}$, i.e., it is close to $1$
and, actually, tends to $1$ as $\alpha\rightarrow0$.

Let now $\alpha\in(\alpha_{0},1)$. In this case $C_{\alpha}=c_{2}(\alpha)$ and
the right-hand side of \eqref{Ball3} is equal to
\[
\frac{2^{\alpha/2}\sqrt{\pi}}{(1-\alpha)(2-\alpha)}\cdot\frac{\Gamma
(\alpha/{2})}{\Gamma((1-\alpha)/{2})}q^{-\alpha/2}.
\]
Therefore, the ratio we are interested in, equals the fraction
\[
\frac{(1-\alpha)(2-\alpha)\Gamma((1-\alpha)/{2})}{2\pi^{(1-\alpha)/{2}}}.
\]
From elementary properties of the $\Gamma$-function it follows that this
expression tends to $1$ as $\alpha\rightarrow1$.

Since both sides of \eqref{Ball2} are larger than the corresponding ones of
\eqref{Ball3}, we can summarize our findings as saying that inequalities
\eqref{Ball2} and \eqref{Ball3} are not comparable, i.e., neither of them is
weaker or stronger than the other.

\section{Non-commutative results of Nazarov-Podkorytov type\label{s4}}

As it was alluded in the Introduction, Nazarov-Podkorytov type inequalities
can be also formulated and proved in the non-commutative setting.

Let ${H}$ be a separable complex Hilbert space. For every compact operator
${A}$ in ${H}$ let ${s}(A)=\{s_{n}({A})\}_{n=1}^{\infty}$ be the sequence of
$s$-numbers of ${A}$ determined by the Schmidt expansion \cite{GK}. For every
$p,q>0$, the class ${{\mathfrak{S}}}^{p,q}$ consists of all compact operators
${A}:\,{H}\rightarrow{H}$ such that
\[
{\Vert A\Vert}_{p,q}:={\Vert{s}(A)\Vert}_{\ell(p,q)}<\infty,
\]
where $\ell(p,q)$, the discrete version of the function space $L(p,q),$
consists of all sequences ${x}=(x_{n})_{n=1}^{\infty}$ of real numbers such
that the quasi-norm
\[
{\Vert{x}\Vert}_{\ell(p,q)}:=\left(  \sum_{n=1}^{\infty}(x_{n}^{\ast}%
)^{q}n^{\frac{q}{p}-1}\right)  ^{\frac{1}{q}}%
\]
is finite. Here, ${x}^{\ast}=(x_{n}^{\ast})_{n=1}^{\infty}$ denotes the
nonincreasing permutation of the sequence $(|x_{n}|)_{n=1}^{\infty}$.

The classes ${{\mathfrak{S}}}^{p,q},$ $p,q>0,$ are two-sided symmetrically
quasi-normed ideals in the space of all bounded operators in ${H,}$ which
sometimes are referred to as Lorentz ideals. The classical Schatten-von
Neumann ideals ${{\mathfrak{S}}}^{p}$ correspond to the case $p=q,$ i.e.,
${{\mathfrak{S}}}^{p}={{\mathfrak{S}}}^{p,p}$.

We introduce a following discrete variant of the condition (\ref{distrib}).
Let $x=(x_{n})_{n=1}^{\infty}$ and $y=(y_{n})_{n=1}^{\infty}$ be two sequences
of real numbers. We shall write $x\lessgtr{y}$ (or $y\gtrless{x}$) if for some
$n_{0}\in\mathbb{N}$ we have
\begin{equation}
x_{n}^{\ast}\leqslant y_{n}^{\ast}\quad\mbox{ if }n\leqslant n_{0}%
\quad\mbox{ and }\quad x_{n}^{\ast}\geqslant y_{n}^{\ast}\quad
\mbox{ if }n>n_{0}. \label{order for seq}%
\end{equation}
Correspondingly, if ${A}$ and ${B}$ are compact operators in ${H}$, we write
${A}\lessgtr{B}$ iff ${s}(A)\lessgtr{s}({B})$.

The main result of this section is a non-commutative version of Theorem
\ref{prop of Dlta}.
For any $p_{0},q_{0}>0$ we denote
\[
\tilde{\Delta}(p_{0},q_{0})=\{(p,q):\,0<p\leqslant p_{0},\,\,0<q\leqslant
q_{0},\,q_{0}/p_{0}\leqslant q/p\}.
\]

\begin{theorem}
\label{IdealTheorem} Let $p_{0},q_{0}>0$, and let ${A}$ and ${B}$ be bounded
compact operators on a complex separable Hilbert space ${H}$ such that
${A}\lessgtr{B}$ and $\Vert{A}\Vert_{p_{0},q_{0}}\geqslant\Vert{B}\Vert
_{p_{0},q_{0}}$ (resp. $\Vert{A}\Vert_{p_{0},q_{0}}>\Vert{B}\Vert_{p_{0}%
,q_{0}}$). Then, the condition $(p,q)\in\tilde{\Delta}(p_{0},q_{0})$ implies
that
\[
{\Vert{A}\Vert}_{p,q}\geqslant{\Vert{B}\Vert}_{p,q}\;(\mbox{resp.}{\Vert
{A}\Vert}_{p,q}>{\Vert{B}\Vert}_{p,q}).
\]

Conversely, if $(p,q)\;{\not \in }\;\tilde{\Delta}(p_{0},q_{0})$, then there
exist two-dimensional operators ${A}$ and ${B}$ such that
\[
{A}\lessgtr{B},\quad{\Vert{A}\Vert}_{p_{0},q_{0}}>{\Vert{B}\Vert}_{p_{0}%
,q_{0}}\quad\mbox{and}\quad{\Vert{A}\Vert}_{p,q}<{\Vert{B}\Vert}_{p,q}.
\]

\end{theorem}

Theorem \ref{IdealTheorem} is an immediate consequence of the above
definitions and the following result for the $\ell(p,q)$-spaces.

\begin{proposition}
\label{SeqTheorem} Let $p_{0},q_{0}>0$.

If ${x}\lessgtr{y}$ and $\Vert{x}\Vert_{\ell(p_{0},q_{0})} \geqslant\Vert
{y}\Vert_{\ell(p_{0},q_{0})}$ (resp. $\Vert{x} \Vert_{\ell(p_{0},q_{0})}%
>\Vert{y}\Vert_{\ell(p_{0},q_{0})}$), then the condition $(p,q)\in
\tilde{\Delta}(p_{0},q_{0})$ implies that
\[
{\Vert{x}\Vert}_{\ell(p,q)}\geqslant{\Vert{y}\Vert}_{\ell(p,q)}%
\;(\mbox{resp.}{\Vert{x}\Vert}_{\ell(p,q)}>{\Vert{y}\Vert}_{\ell(p,q)}).
\]

Conversely, if $(p,q)\;{\not \in }\;\tilde{\Delta}(p_{0},q_{0})$, then there
exist sequences $x$ and $y$, having at most two non-zero entries, such that
\[
{x}\lessgtr{y},\quad{\Vert{x}\Vert}_{\ell(p_{0},q_{0})}>{\Vert{y}\Vert}%
_{\ell(p_{0},q_{0})}\quad\mbox{and}\quad{\Vert{x}\Vert}_{\ell(p,q)}<{\Vert
{y}\Vert}_{\ell(p,q)}.
\]

\end{proposition}

\begin{proof}
We consider first the case when $q/p=q_{0}/p_{0}$ and $q<q_{0}$. Let
$r=q/q_{0}$ and let $u_{0}$ be the solution of the equation $f^{\prime}(u)=1$,
where $f(u)=u^{r}$, i.e., $u_{0}={r}^{\frac{1}{1-r}}$. Therefore, since $r<1$,
$f^{\prime}(u)$ decreases. It follows that for $v\geqslant u\geqslant u_{0}$
we have
\begin{equation}
v^{r}-u^{r}\leqslant v-u, \label{ineq1 for seq}%
\end{equation}
while for $v\leqslant u\leqslant u_{0},$ it holds
\begin{equation}
v^{r}-u^{r}\geqslant v-u. \label{ineq2 for seq}%
\end{equation}

Let now ${x}\lessgtr{y}$ and $\Vert{x}\Vert_{\ell(p_{0},q_{0})} \geqslant
\Vert{y}\Vert_{\ell(p_{0},q_{0})}$. Without loss of generality, we can assume
that ${x}={x}^{\ast}$, ${y}={y}^{\ast}$. Moreover, observe that for each $s>0$
the inequalities $\Vert s{x}\Vert_{\ell(p,q)}\leq\Vert s{y}\Vert_{\ell(p,q)}$
and $\Vert{x}\Vert_{\ell(p,q)}\leq\Vert{y}\Vert_{\ell(p,q)}$ (resp. $\Vert
s{x}\Vert_{\ell(p,q)}<\Vert s{y}\Vert_{\ell(p,q)}$ and $\Vert{x}\Vert
_{\ell(p,q)}<\Vert{y}\Vert_{\ell(p,q)}$) are equivalent. Therefore, we can
also assume that $x_{n_{0}}^{q_{0}}=u_{0}$. Then, using inequalities
\eqref{order for seq}, \eqref{ineq1 for seq} and \eqref{ineq2 for seq}, we
obtain
\begin{align*}
{\Vert{x}\Vert}_{\ell(p,q)}^{q}-{\Vert{y}\Vert}_{\ell(p,q)}^{q}  &
=\sum_{n=n_{0}+1}^{\infty}(x_{n}^{q}-y_{n}^{q})n^{\frac{q}{p}-1}-\sum
_{n=1}^{n_{0}}(y_{n}^{q}-x_{n}^{q})n^{\frac{q}{p}-1}\\
&  =\sum_{n=n_{0}+1}^{\infty}\bigl((x_{n}^{q_{0}})^{r}-(y_{n}^{q_{0}}%
)^{r}\bigr)n^{\frac{q}{p}-1}-\sum_{n=1}^{n_{0}}\bigl((y_{n}^{q_{0}}%
)^{r}-(x_{n}^{q_{0}})^{r}\bigr)n^{\frac{q}{p}-1}\\
&  \geqslant\sum_{n=n_{0}+1}^{\infty}(x_{n}^{q_{0}}-y_{n}^{q_{0}}%
)n^{\frac{q_{0}}{p_{0}}-1}-\sum_{n=1}^{n_{0}}(y_{n}^{q_{0}}-x_{n}^{q_{0}%
})n^{\frac{q_{0}}{p_{0}}-1}\\
&  ={\Vert{{x}\Vert}}_{\ell(p_{0},q_{0})}^{q_{0}}-{\Vert{{y}\Vert}}%
_{\ell(p_{0},q_{0})}^{q_{0}}\geqslant0.
\end{align*}

Consider now the case $q=q_{0}$, $p<p_{0}$. Since $\frac{q}{p}-\frac{q_{0}%
}{p_{0}}>0$, we have
\[
n^{\frac{q}{p}-\frac{q_{0}}{p_{0}}}\leqslant n_{0}^{\frac{q}{p}-\frac{q_{0}%
}{p_{0}}}\quad\mbox{ if }n\leqslant n_{0},\quad\mbox{ and }\quad n^{\frac
{q}{p}-\frac{q_{0}}{p_{0}}}\geqslant n_{0}^{\frac{q}{p}-\frac{q_{0}}{p_{0}}%
}\quad\mbox{ if }n>n_{0}.
\]
Hence,
\begin{align*}
{\Vert{{x}\Vert}}_{\ell(p,q)}^{q}-{\Vert{{y}\Vert}}_{\ell(p,q)}^{q}  &
=\sum_{n=n_{0}+1}^{\infty}(x_{n}^{q}-y_{n}^{q})n^{\frac{q}{p}-1}-\sum
_{n=1}^{n_{0}}(y_{n}^{q}-x_{n}^{q})n^{\frac{q}{p}-1}\\
&  =\sum_{n=n_{0}+1}^{\infty}(x_{n}^{q_{0}}-y_{n}^{q_{0}})n^{\frac{q_{0}%
}{p_{0}}-1}\cdot n^{\frac{q}{p}-\frac{q_{0}}{p_{0}}}-\sum_{n=1}^{n_{0}}%
(y_{n}^{q_{0}}-x_{n}^{q_{0}})n^{\frac{q_{0}}{p_{0}}-1}\cdot n^{\frac{q}%
{p}-\frac{q_{0}}{p_{0}}}\\
&  \geqslant\sum_{n=n_{0}+1}^{\infty}(x_{n}^{q_{0}}-y_{n}^{q_{0}}%
)n^{\frac{q_{0}}{p_{0}}-1}\cdot n_{0}^{\frac{q}{p}-\frac{q_{0}}{p_{0}}}%
-\sum_{n=1}^{n_{0}}(y_{n}^{q_{0}}-x_{n}^{q_{0}})n^{\frac{q_{0}}{p_{0}}-1}\cdot
n_{0}^{\frac{q}{p}-\frac{q_{0}}{p_{0}}}\\
&  =n_{0}^{\frac{q}{p}-\frac{q_{0}}{p_{0}}}\bigl({\Vert{{x}\Vert}}_{\ell
(p_{0},q_{0})}^{q_{0}}-{\Vert{{y}\Vert}}_{\ell(p_{0},q_{0})}^{q_{0}%
}\bigr)\geqslant0.
\end{align*}

A similar argument shows that ${\Vert{{x}\Vert}}_{\ell(p_{0},q_{0})}%
>{\Vert{{y}\Vert}}_{\ell(p_{0},q_{0})}$ implies ${\Vert{{x}\Vert}}_{\ell
(p,q)}>{\Vert{{y}\Vert}}_{\ell(p,q)}$.

Finally, for any $(p,q)\in\tilde{\Delta}(p_{0},q_{0})$, the desired result can
be achieved by a transitivity argument (see Lemma \ref{lemma trans1}), i.e.,
passing successively from $(p_{0},q_{0})$ to $(p_{0}q/q_{0},q)$ and then to
$(p,q)$.

To prove the second assertion, suppose that $(p,q)\;{\not \in }\;\tilde
{\Delta}(p_{0},q_{0})$. Then, $(p,q)\in F_{1}\bigcup F_{2}\bigcup F_{3},$
where the sets $F_{1}$, $F_{2}$ and $F_{3}$ are defined by
\[
F_{1}=\{(p,q):p>0,\,q>q_{0}\},\quad\quad F_{2}=\{(p,q):\,q=q_{0}%
,\,p>0,\,q/p<q_{0}/p_{0}\},
\]%
\[
F_{3}=\{(p,q):\,p>0,\,0<q<q_{0},q/p<q_{0}/p_{0}\}.
\]

It will be convenient to identify a two-dimensional vector $\vec{x}%
=(x_{1},x_{2})$ with the sequence $x=(x_{n})_{n=1}^{\infty}$ with $x_{n}=0$
for $n>2$. Let
\[
r:=\frac{q}{q_{0}},\quad\gamma:=\frac{q}{p}-1,\quad\gamma_{0}:=\frac{q_{0}%
}{p_{0}}-1,
\]
and consider a two-dimensional vector $\vec{z}(\alpha)=(z_{1}(\alpha
),z_{2}(\alpha))$ with coordinates $z_{i}=z_{i}(\alpha)$, $i=1,2$, selected
such that
\[
z_{1}^{q_{0}}=1-\alpha\cdot2^{\gamma_{0}}\geqslant\alpha
\;\;\mbox{and}\;\;z_{2}^{q_{0}}=\alpha\geqslant0.
\]
Then,
\[
z_{1}\geqslant z_{2}\geqslant0,\quad{\Vert{z}(\alpha)\Vert}_{\ell(p_{0}%
,q_{0})}^{q_{0}}=z_{1}^{q_{0}}+z_{2}^{q_{0}}\cdot2^{\gamma_{0}}=1.
\]
Let
\[
\varphi(\alpha):={\Vert{z}(\alpha)\Vert}_{\ell(p,q)}^{q}=z_{1}^{q}+z_{2}%
^{q}\cdot2^{\gamma}=(1-\alpha\cdot2^{\gamma_{0}})^{r}+\alpha^{r}\cdot
2^{\gamma},
\]
and by computation
\begin{equation}
\varphi^{\prime}(\alpha)=-r2^{\gamma_{0}}(1-\alpha\cdot2^{\gamma_{0}}%
)^{r-1}+r2^{\gamma}\alpha^{r-1}=r2^{\gamma}\alpha^{r-1}\left(  1-2^{\gamma
_{0}-\gamma}\left(  \frac{1}{\alpha}-2^{\gamma_{0}}\right)  ^{r-1}\right)  .
\label{formula}%
\end{equation}

Suppose first that $(p,q)\in F_{1}\cup F_{2}$. Then, either $r>1$, or $r=1$
and $\gamma_{0}-\gamma>0.$ In either of these cases, using (\ref{formula}) we
see that there exists $\delta>0$ such that $\varphi^{\prime}(\alpha)<0$ if
$\alpha\in(0,\delta)$. Therefore, for each $\alpha\in(0,\delta)$ we have
\[
{\Vert{z}(\alpha)\Vert}_{\ell(p,q)}^{q}<{\Vert{z}(0)\Vert}_{\ell(p,q)}^{q}=1.
\]

Next, fix any $\alpha_{0}\in(0,\delta)$ such that $1-\alpha_{0}\cdot
2^{\gamma_{0}}>\alpha_{0}$ or equivalently $\alpha_{0}<{1}/(1+2^{\gamma_{0}}%
)$. For $\varepsilon>0,$ let $x_{1}$ and $x_{2}$ be defined by
\[
x_{1}^{q_{0}}=1-\alpha_{0}\cdot2^{\gamma_{0}},\quad x_{2}^{q_{0}}=\alpha
_{0}+\varepsilon,
\]
Therefore, if we choose $\varepsilon>0$ sufficiently small, then the preceding
discussion shows that the sequences $x$ and $y$ corresponding to the vectors
$\vec{x}=(x_{1},x_{2})$ and $\vec{y}=(1,0)$ satisfy the required inequalities
\[
{{x}\lessgtr{y},\quad{\Vert{x}\Vert}_{\ell(p_{0},q_{0})}>{\Vert{y}\Vert}%
_{\ell(p_{0},q_{0})}\quad\mbox{and}\quad{\Vert{x}\Vert}_{\ell(p,q)}<{\Vert
{y}\Vert}_{\ell(p,q)}.}%
\]

Assume now that $(p,q)\in F_{3}$. Then, we have
\begin{equation}
r<1\quad\mbox{and}\quad\gamma_{0}-\gamma>0. \label{IneqF1}%
\end{equation}
Using (\ref{formula}) and a direct calculation, we find that the root
$\alpha_{1}$ of the equation $\varphi^{\prime}(\alpha)=0$, is given by
\[
\alpha_{1}=\frac{1}{2^{\frac{\gamma_{0}-\gamma}{1-r}}+2^{\gamma_{0}}}.
\]
Moreover, from \eqref{IneqF1} it follows that $\varphi$ attains a strict
maximum at $\alpha_{1}$ and $\frac{\gamma_{0}-\gamma}{1-r}>0$. Consequently,
\[
1-\alpha_{1}\cdot2^{\gamma_{0}}=1-\frac{2^{\gamma_{0}}}{2^{\frac{\gamma
_{0}-\gamma}{1-r}}+2^{\gamma_{0}}}=\frac{2^{\frac{\gamma_{0}-\gamma}{1-r}}%
}{2^{\frac{\gamma_{0}-\gamma}{1-r}}+2^{\gamma_{0}}}=2^{\frac{\gamma_{0}%
-\gamma}{1-r}}\alpha_{1}>\alpha_{1}.
\]
Furthermore, if $\alpha_{2}:=1/(1+2^{\gamma_{0}})$, then $\alpha_{1}%
<\alpha_{2}<1$. Thus, setting
\[
x_{1}^{q_{0}}=x_{2}^{q_{0}}=\alpha_{2}+\varepsilon,\quad y_{1}^{q_{0}%
}=1-\alpha_{1}\cdot2^{\gamma_{0}}\quad\mbox{and}\quad y_{2}^{q_{0}}=\alpha
_{1},
\]
we see that selecting $\varepsilon>0$ sufficiently small, the sequences $x$
and $y$, corresponding to the vectors $\vec{x}=(x_{1},x_{2})$ and $\vec
{y}=(y_{1},y_{2})$, satisfy the desired inequalities
\[
{{x}\lessgtr{y},\quad{\Vert{x}\Vert}_{\ell(p_{0},q_{0})}>{\Vert{y}\Vert}%
_{\ell(p_{0},q_{0})}\quad\mbox{and}\quad{\Vert{x}\Vert}_{\ell(p,q)}<{\Vert
{y}\Vert}_{\ell(p,q)}.}%
\]
This completes the proof of the proposition.
\end{proof}

\begin{remark}
A result similar to Proposition \ref{SeqTheorem} holds also for the
$N$-dimensional spaces $\ell_{N}(p,q)$, $N\geqslant2$, equipped with the
quasi-norms
\[
{\Vert{x}\Vert}_{\ell_{N}(p,q)}:=\left(  \sum_{n=1}^{N}(x_{n}^{\ast}%
)^{q}n^{\frac{q}{p}-1}\right)  ^{\frac{1}{q}}.
\]

\end{remark}

Precisely in the same way we can obtain the following counterpart of Theorem
\ref{IdealTheorem} (cf. Theorem \ref{prop of Dlta}).

\begin{theorem}
\label{IdealTheorem1} Let $p_{0},q_{0}>0$, and let ${A}$ and ${B}$ be bounded
compact operators on a complex separable Hilbert space ${H}$ such that
${A}\gtrless{B}$ and $\Vert{A}\Vert_{p_{0},q_{0}}\geqslant\Vert{B}\Vert
_{p_{0},q_{0}}$ (resp. $\Vert{A}\Vert_{p_{0},q_{0}}>\Vert{B}\Vert_{p_{0}%
,q_{0}}$). Let
\[
{\Delta}(p_{0},q_{0}):=\{(p,q):\,p\geqslant p_{0},\,q\geqslant q_{0}%
,\,q_{0}/p_{0}\geqslant q/p\}.
\]
Then, the condition $(p,q)\in{\Delta}(p_{0},q_{0})$ implies that
\[
{\Vert{A}\Vert}_{p,q}\geqslant{\Vert{B}\Vert}_{p,q}\;(\mbox{resp.}{\Vert
{A}\Vert}_{p,q}>{\Vert{B}\Vert}_{p,q}).
\]

Conversely, if $(p,q)\;{\not \in }\;{\Delta}(p_{0},q_{0})$, then there exist
two-dimensional operators ${A}$ and ${B}$ such that
\[
{A}\gtrless{B},\quad{\Vert{A}\Vert}_{p_{0},q_{0}}>{\Vert{B}\Vert}_{p_{0}%
,q_{0}}\quad\mbox{and}\quad{\Vert{A}\Vert}_{p,q}<{\Vert{B}\Vert}_{p,q}.
\]

\end{theorem}

\section{Interpolation methods\label{Int meth}}

For further background information we refer the reader to \cite{BL},
\cite{BS}, \cite{LT} and \cite{KPS}.

\subsection{The real $K-$method\label{K-method}}

As we pointed out in the Introduction (cf. \eqref{intro2}), one of the
possible equivalent formulations of the Hardy-Littlewood-Polya order can be
written using the maximal averages of $f$ and $g$, i.e., the functions
\[
\int_{0}^{t}f^{\ast}(s)\,ds=K(t,f;L^{1},L^{\infty})\;\;\mbox{and}\;\;\int%
_{0}^{t}g^{\ast}(s)\,ds=K(t,g;L^{1},L^{\infty}).
\]
More generally, we will consider majorization conditions associated to the
$K-$functional of a given pair $(X,Y)$ of compatible\footnote{We shall say
that a pair $(X,Y)$ of Banach (or quasi-Banach) spaces is compatible if $X$
and $Y$ are both linearly and continuously embedded in some Hausdorff linear
topological space. For example, if $(X,Y)$ is an \textquotedblleft ordered
pair" (that is, e.g. $Y\subset X$), then the pair is obviously compatible.}
Banach (quasi-Banach) spaces. We define the Peetre $K-$functional for $h\in
X+Y,$ $t>0,$ by
\[
K(t,h;X,Y)=\inf\{\left\Vert f_{1}\right\Vert _{X}+t\left\Vert f_{2}\right\Vert
_{Y}:h=f_{1}+f_{2},f_{1}\in X,f_{2}\in Y\},\text{ \ \ }t>0.
\]
The $K-$functional itself provides an interpolation between $X$ and $Y$ that
corresponds to the computations%
\[
\lim_{t\rightarrow\infty}K(t,h;X,Y)=\left\Vert h\right\Vert _{X}%
,\;\;\;\lim_{t\rightarrow0}\frac{K(t,h;X,Y)}{t}=\left\Vert h\right\Vert _{Y},
\]
which are valid whenever the spaces $X$ and $Y$ are Gagliardo complete with
respect to the sum $X+Y$ (see e.g. \cite[p.~295-296]{BS}). This leads to the
corresponding extended majorization of the Hardy-Littlewood-Polya type
$\succeq_{K(X,Y)}$ defined by%
\[
f\succeq_{K(X,Y)}g\Leftrightarrow\text{ }K(t,f;X,Y)\geq K(t,g;X,Y)\text{ for
all }t>0,
\]
where the classical case corresponds to the pair $(X,Y)=(L^{1},L^{\infty}).$
Likewise, for $t_{1}>0,$ we define%
\[
f\succeq_{t_{1}K(X,Y)}g\text{ }\Leftrightarrow\text{ }\left\{
\begin{array}
[c]{cc}%
K(t,f;X,Y)\geq & K(t,g;X,Y)\text{ for }t<t_{1},\\
K(t,f;X,Y)\leq & K(t,g;X,Y)\text{ for }t\geq t_{1}.
\end{array}
\right.
\]

Recall that the Lions-Peetre interpolation spaces $(X,Y)_{\theta,q}$,
$0<\theta<1$, $0<q\leq\infty,$ are defined by
\[
(X,Y)_{\theta,q}=\{f:\left\Vert f\right\Vert _{(X,Y)_{\theta,q}}:=\left\{
\int_{0}^{\infty}[t^{-\theta}K(t,f;X,Y)]^{q}\frac{dt}{t}\right\}
^{1/q}<\infty\},\text{ }%
\]
with the usual modification when $q=\infty.$

The natural corresponding \textquotedblleft extra\textquotedblright%
\ cancellation condition in this case corresponds to demand that for some
$(\theta_{0},q_{0})\in(0,1)\times(0,\infty]$ we have%
\[
\left\Vert f\right\Vert _{(X,Y)_{\theta_{0},q_{0}}}\geq\left\Vert g\right\Vert
_{(X,Y)_{\theta_{0},q_{0}}}.
\]

In what follows when dealing with the $K-$functional of an element $f$ and
when the pair $(X,Y)$ is understood, we will simplify the notation and just
write $K(t,f).$

\begin{theorem}
\label{teoderivado} Let $0<\theta,\eta<1$, $0<p\leq q<\infty,$ satisfy the
inequality $q(1-\eta)\leq p(1-\theta)$. Let $(X,Y)$ be a Banach pair, $f,g\in
X+Y$. Suppose that there exists $t_{1}>0$ such that%
\begin{equation}
f\succeq_{t_{1}K(X,Y)}g, \label{eq101}%
\end{equation}

Then, from
\begin{equation}
\left\Vert g\right\Vert _{(X,Y)_{\theta,p}}\leq\left\Vert f\right\Vert
_{(X,Y)_{\theta,p}} \label{eq103}%
\end{equation}
it follows that
\begin{equation}
\left\Vert g\right\Vert _{(X,Y)_{\eta,q}}\leq\left\Vert f\right\Vert
_{(X,Y)_{\eta,q}}. \label{condo3'}%
\end{equation}

\end{theorem}

\begin{proof}
Suppose first that $q=p$. Then, $0<\theta\leq\eta$.

Let
\[
\phi_{f,g}(\eta):=t_{1}^{\eta p}\left\Vert f\right\Vert _{(X,Y)_{\eta,p}}%
^{p}-t_{1}^{\eta p}\left\Vert g\right\Vert _{(X,Y)_{\eta,p}}^{p},
\]
and compute%
\begin{align*}
\phi_{f,g}(\eta)-\phi_{f,g}(\theta)  &  =\int_{0}^{\infty}\left(
K(s,f)^{p}-K(s,g)^{p}\right)  \left(  \frac{s}{t_{1}}\right)  ^{-p\eta}%
\,\frac{ds}{s}\\
&  -\int_{0}^{\infty}\left(  K(s,f)^{p}-K(s,g)^{p}\right)  \left(  \frac
{s}{t_{1}}\right)  ^{-p\theta}\,\frac{ds}{s}\\
&  =\int_{0}^{\infty}\left(  K(s,f)^{p}-K(s,g)^{p}\right)  \left(  \left(
\frac{s}{t_{1}}\right)  ^{-p\eta}-\left(  \frac{s}{t_{1}}\right)  ^{-p\theta
}\right)  \,\frac{ds}{s}.
\end{align*}
Suppose that $s\in(0,t_{1}).$ Then, $({s}/{t_{1}})^{-p\eta}$ $\ $\ is an
increasing function of $\eta$, and $({s}/{t_{1}})^{-p\eta}-({s}/{t_{1}%
})^{-p\theta}\geq0$ for $\eta\geq\theta.$ This fact combined with assumption
\eqref{eq101} yields that for $s\in(0,t_{1})$ we have
\begin{equation}
\left(  K(s,f)^{p}-K(s,g)^{p}\right)  \cdot\left(  ({s}/{t_{1}})^{-p\eta}%
-({s}/{t_{1}})^{-p\theta}\right)  \geq0. \label{signo}%
\end{equation}

Similarly, if $s\in(t_{1},\infty)$ then $({s}/{t_{1}})^{-p\eta}$ $\ $\ is a
decreasing function of $\eta,$ and (\ref{signo}) holds since now both factors
are negative. Thus, it follows that%
\[
\phi_{f,g}(\eta)-\phi_{f,g}(\theta)\geq0.
\]
In addition, by \eqref{eq103}, we have $\phi_{f,g}(\theta)\geq0$. Combining
inequalities, we obtain that $\phi_{f,g}(\eta)\geq0$, which is equivalent to
desired inequality \eqref{condo3'}. This concludes the proof in the case when
$q=p$.

Let now $q>p$. We will show that matters can be reduced to an application of
Proposition \ref{cor1b} (a).

One can easily check that $\eta+\alpha/q=\theta+\alpha/p$ for
\[
\alpha=\frac{(\eta-\theta)pq}{q-p}.
\]
Moreover, from the inequality $q(1-\eta)\leq p(1-\theta)$ it follows that
$\alpha\geq p(1-\theta),$ and hence $\theta+\alpha/p\geq1$.

Let $h\in X+Y$. Define the function $u_{h}(t),$ for $t>0,$ by
\[
u_{h}(t):=\frac{K(t,h)}{t^{\theta+\alpha/p}}=\frac{K(t,h)}{t^{\eta+\alpha/q}%
}.
\]
Since $K(t,h)/t$ is a nonincreasing function in $t$ and $\eta+\alpha
/q=\theta+\alpha/p\geq1$, it follows that $u_{h}(t)$ is a decreasing function.
By the definition of the Lions-Peetre interpolation spaces and the Lorentz
$L(p,q)$-spaces (cf. (\ref{quasi-norm*})), we can write%
\begin{align*}
\left\Vert h\right\Vert _{(X,Y)_{\theta,p}}^{p}  &  =\int_{0}^{\infty}\left(
\frac{K(t,h)}{t^{\theta+\alpha/p}}\right)  ^{p}t^{\alpha-1}\,dt\\
&  =\int_{0}^{\infty}\left(  u_{h}(t)\right)  ^{p}t^{\alpha-1}\,dt\\
&  =\Vert u_{h}\Vert_{L(p/\alpha,p)}^{p},
\end{align*}
and likewise,%
\[
\left\Vert h\right\Vert _{(X,Y)_{\eta,q}}^{q}=\Vert u_{h}\Vert_{L(q/\alpha
,q)}^{q}.
\]
Moreover, hypothesis \eqref{eq101} yields
\[
u_{g}(t)\leq u_{f}(t)\text{ if }t\leq t_{1},\;\;\mbox{and}\;\;u_{g}(t)\geq
u_{f}(t)\;\;\mbox{ if }\;t>t_{1}.
\]
Therefore, since $u_{f}(t)$ and $u_{g}(t)$ are decreasing functions, there
exists $\tau_{0}>0$ such that
\[
\lambda_{u_{f}}(\tau)\left\{
\begin{array}
[c]{cc}%
\leq & \lambda_{u_{g}}(\tau)\text{ for }\tau<\tau_{0}\\
\geq & \lambda_{u_{g}}(\tau)\text{ for }\tau>\tau_{0}.
\end{array}
\right.
\]
Moreover, from \eqref{eq103} we have $\Vert u_{f}\Vert_{L(p/\alpha,p)}%
\geq\Vert u_{g}\Vert_{L(p/\alpha,p)}$. Since $q>p$, we can apply Proposition
\ref{cor1b}(a), to obtain
\[
\left\Vert f\right\Vert _{(X,Y)_{\eta,q}}^{q}=\Vert u_{f}\Vert_{L(q/\alpha
,q)}^{q}\geq\Vert u_{g}\Vert_{L(q/\alpha,q)}^{q}=\left\Vert g\right\Vert
_{(X,Y)_{\eta,q}}^{q},
\]
and the desired result follows.
\end{proof}

In particular, if $(X,Y)$ is the pair $(L^{1},L^{\infty})$ we get

\begin{corollary}
Let $0<\theta,\eta<1$, $0<p\leq q<\infty,$ be so that $q(1-\eta)\leq
p(1-\theta)$. Suppose that functions $f,g\in L^{1}+L^{\infty}$ satisfy:
\[
\int_{0}^{t} f^{*}(s)\,ds\left\{
\begin{array}
[c]{cc}%
\geq & \int_{0}^{t} g^{*}(s)\,ds\text{ for }t<t_{1},\\
\leq & \int_{0}^{t} g^{*}(s)\,ds\text{ for }t>t_{1}%
\end{array}
\right.
\]
for some $t_{1}>0$, and
\[
\left\Vert g\right\Vert _{(L^{1},L^{\infty})_{\theta,p}}\leq\left\Vert
f\right\Vert _{(L^{1},L^{\infty})_{\theta,p}}.
\]
Then,
\[
\left\Vert g\right\Vert _{(L^{1},L^{\infty})_{\eta,q}}\leq\left\Vert
f\right\Vert _{(L^{1},L^{\infty})_{\eta,q}}.
\]

\end{corollary}

\begin{remark}
Suppose $0<\eta<1$, $0<q<\infty$ are arbitrary. It is well known (see e.g.
\cite[Theorem~5.2.1]{BL}) that $(L^{1},L^{\infty})_{\eta,q}=L(p,q)$ (with
equivalence of norms), where $p=1/(1-\eta)$. Therefore, one can see that the
result of the last corollary differs sharply from that of Theorem
\ref{prop of Dlta}, where the spaces $L(p,q)$ are defined using the
quasi-norms \eqref{quasi-norm}.
\end{remark}

\begin{remark}
\label{modified} Let $a>0$. Similar results hold also for the modified
interpolation spaces\footnote{See Section \ref{derivatives of norms}.}
$\langle X,Y\rangle_{\theta,q}$ defined by
\[
\langle X,Y\rangle_{\theta,q}=\{f:\left\Vert f\right\Vert _{\langle
X,Y\rangle_{\theta,q}}:=\left\{  \int_{0}^{a}[t^{-\theta}K(t,f;X,Y)]^{q}%
\frac{dt}{t}\right\}  ^{1/q}<\infty\}.
\]

\end{remark}

\subsection{The $E$-method}

Similarly we can formulate results of Nazarov-Podkorytov type in the context
of the $E-$method of interpolation. In this section we give only a brief
outline of this direction of research.

Recall that the distribution function and the non-increasing rearrangement are
special instances of $E$-functionals for suitable pairs of function spaces.
Indeed, let
\begin{equation}
E(t,f;L^{\infty},L^{0})=\inf_{\left\Vert h\right\Vert _{L^{\infty}}\leq
t}\{\left\Vert f-h\right\Vert _{L^{0}}\}, \label{e1}%
\end{equation}
where $L^{0}$ denotes the set of all measurable functions $f$ with finite
support, i.e., such that%
\begin{equation}
\left\Vert f\right\Vert _{L^{0}}=\mu{\{x:\,f(x)\neq0\}}<\infty. \label{e1a}%
\end{equation}
Then, by a simple computation (cf. \cite{BL}, \cite{pee0}), we find
\[
E(t,f;L^{\infty},L^{0})=\lambda_{f}(t).
\]
Likewise, reversing the order of the spaces, we obtain the generalized inverse
of $\lambda_{f},$ i.e. $f^{\ast}:$%
\[
E(t,f;L^{0},L^{\infty})=\inf_{\left\Vert h\right\Vert _{L^{0}}\leq
t}\{\left\Vert f-h\right\Vert _{L^{\infty}}\}=f^{\ast}(t).
\]
At this point we see that the basic elements of the Nazarov-Podkorytov method
for $L(p,q)$-spaces can be easily extended to general approximation spaces as
follows. For a compatible pair of quasi-Banach spaces $(X,Y)$ we define the
\textquotedblleft error\textquotedblright\ functional $E(t,f;X,Y)$ by (cf.
\cite[Chapter~7]{BL})%
\[
E(t,f;X,Y):=\inf_{\left\Vert h\right\Vert _{X}\leq t}\{\left\Vert
f-h\right\Vert _{Y}\},
\]
and the corresponding spaces using the quasi-norms given by (cf.
\cite[Chapter~7]{BL}, \cite{pee0})%
\[
\left\Vert f\right\Vert _{(X,Y)_{p,q;E}}=\left\{  \int_{0}^{\infty}\left(
E(t,f;X,Y)\right)  ^{q/p}d(t^{q})\right\}  ^{1/q},\;\;0<p,q<\infty.
\]

The Nazarov-Podkorytov type results for $(X,Y)_{p,q;E}$ spaces now follow
\textit{mutatis mutandis} by means of reformulating the assumption
(\ref{distrib}) in terms of $E-$functionals. This approach subsumes the
results presented in the previous sections concerning $L(p,q)$-spaces and
their non-commutative counterparts. We shall leave the details to the
interested readers.

\section{Going beyond $L^{p}$ inequalities: Majorization via $K-$functionals
\label{s7}}

As we have seen, $K-$functionals can be used to provide a far reaching
extension of the classical theory of majorization. Here, we use further ideas
from interpolation theory to establish Nazarov-Podkorytov type results using
majorization with respect to some natural functionals appearing in this theory.

We have observed that the connection between majorization and interpolation is
provided by the \textit{exact} formula for the $K-$functional for the pair
$(L^{1},L^{\infty})$ (see \eqref{K-f for}).
However, $K-$functionals are generally known only up to constants. For
example, the corresponding formula for the pair $(L^{p_{0}},L^{\infty})$,
$p_{0}>1$, takes the form (cf. \cite[Theorem~5.2.1]{BL})%
\begin{equation}
K(t,f;L^{p_{0}},L^{\infty})\asymp\left(  \int_{0}^{t^{p_{0}}}f^{\ast
}(s)^{p_{0}}\,ds\right)  ^{1/p_{0}},\;\;t>0, \label{K-functLp}%
\end{equation}
where $\asymp$ indicates the existence of two-sided estimates with constants.
This leads to the following definition:%
\[
f\succeq_{p_{0}}g\Longleftrightarrow\int_{0}^{t}f^{\ast}(s)^{p_{0}}%
\,ds\geq\int_{0}^{t}g^{\ast}(s)^{p_{0}}\,ds,\;\;t>0.
\]

Let $f,g\in L^{p_{0}}$ satisfy condition \eqref{distrib} and suppose moreover
that inequality ${\Vert f\Vert}_{p_{0}}\geq{\Vert g\Vert}_{p_{0}}$ holds.
Then, by Lemma \ref{lemma:informal}, we find that
\begin{equation}
\int\limits_{0}^{t}f^{\ast}(s)^{p_{0}}\,ds\geq\int\limits_{0}^{t}g^{\ast
}(s)^{p_{0}}\,ds\;\;\mbox{for all}\;\;t>0. \label{IntIneq-a}%
\end{equation}
Therefore, if $X$ is a rearrangement invariant space endowed with a norm of
the form
\begin{equation}
{\Vert h\Vert}_{X}=\biggl\|\Bigl(\int_{0}^{t^{p_{0}}}h^{\ast}(s)^{p_{0}%
}\,ds\Bigr)^{1/p_{0}}\biggr\|_{F}, \label{IntLp}%
\end{equation}
where $F$ is a Banach lattice on $(0,\infty)$, it follows that
(\ref{IntIneq-a}) implies%
\begin{equation}
\Vert f\Vert_{X}\geq\Vert g\Vert_{X} \label{deduce}%
\end{equation}
(with the convention that $\Vert h\Vert_{X}=\infty$ if $h\not \in X$).

The spaces $X$ that can be normed by expressions of the form (\ref{IntLp}) are
precisely characterized using interpolation theory. Indeed, since $(L^{p_{0}
},L^{\infty})$ is a $K$-monotone or Calder\'{o}n-Mityagin pair \cite{LSh} (see
also \cite{sparr}), by a well-known result due to Brudnyi-Kruglyak
\cite[Theorem~4.4.5]{BK91}, all interpolation spaces between $L^{p_{0}}$ and
$L^{\infty}$ admit a renorming of the form \eqref{IntLp}. This application of
abstract interpolation gives an extension of the Nazarov-Podkorytov Lemma to
interpolation spaces for the pair $(L^{p_{0}},L^{\infty})$, modulo an
appropriate renorming. We now show that for a suitable class of interpolation
spaces renorming is not necessary.

Let $1<p<\infty.$ A Banach lattice $E$ is said to be \textit{$p$-convex}, if
there exists a constant $C>0$ such that, for all sequences $\{f_{k}%
\}_{k=1}^{n}\subset E,$ we have
\begin{equation}
\biggl\|\Bigl(\sum_{k=1}^{n}|f_{k}|^{p}\Bigr)^{1/p}\biggr\|_{E}\leq
C\Bigl(\sum_{k=1}^{n}{\Vert f_{k}\Vert}_{E}^{p}\Bigr)^{1/p}. \label{arriba}%
\end{equation}
If the constant $C$ can be chosen to be $1,$ then we say that $E$ is
\textit{exact $p$-convex. }It is well known that every $p$-convex
rearrangement invariant space is an interpolation space with respect to the
pair $(L^{p},L^{\infty})$, $1<p<\infty$ (cf. \cite[Theorem~7.3]{KM-Handbook}).
Moreover,
if $X$ is an exact $p$-convex rearrangement invariant space, then we can write
${\Vert f\Vert}_{X}={\Vert|f|^{p}\Vert}_{Y}^{1/p}$, where $Y$ is a
rearrangement invariant space (see e.g. \cite[1.d]{LT}). Therefore, if $X$ is
exact $p_{0}-$convex and \eqref{IntIneq-a} holds for some functions $f\in X$
and $g$, we achieve estimate (\ref{deduce}) without renormings. In particular,
this result contains the original Nazarov-Podkorytov Lemma since one can
easily verify that $L^{q}$ is an exact $p$-convex rearrangement invariant
space if and only if $q\geq p$.

Summarizing, we have obtained the following result.

\begin{theorem}
\label{cor1c} Let $1\leq p<\infty$.

(i) Suppose $X$ is an interpolation space between $L^{p}$ and $L_{\infty}$.
Then $X$ admits a renorming $\Vert\cdot\Vert_{X}^{\prime}$ such that, if
functions $f\in L^{p}$ and $g$ satisfy condition \eqref{distrib} and ${\Vert
f\Vert}_{p}\geq{\Vert g\Vert}_{p}$, it follows that $\Vert f\Vert_{X}^{\prime
}\geq\Vert g\Vert_{X}^{\prime}$.

(ii) If, moreover, $X$ is an exact $p$-convex rearrangement invariant space,
then under the same assumptions we actually get $\Vert f\Vert_{X}\geq\Vert
g\Vert_{X}$.
\end{theorem}

In particular, Theorem \ref{cor1c}(ii) proves the extension of Ball's integral
inequality (\ref{intro1}) stated in Theorem \ref{teoexact}. Likewise, it is
well known (see e.g. \cite{KMP97}) that an Orlicz space $L_{M}$ equipped with
the Luxemburg norm (see \eqref{luxemburg}) is exact $p$-convex if and only if
the function $M({u}^{1/p})$, $u>0$, is convex. This observation combined with
Theorem \ref{teoexact} proves Corollary \ref{teoexact:cor}.

\begin{remark}
It is important to note here that the class of rearrangement invariant spaces
that are interpolation spaces between $L^{p_{0}}$ and $L^{\infty}$ is much
wider than the class of exact $p_{0}$-convex rearrangement invariant spaces.
For example, the Lorentz space $L(p,1)$, with $p>p_{0}$, is an interpolation
space between $L^{p_{0}}$ and $L^{\infty}$, that is not $r$-convex for any
$r>1$ (cf. e.g. \cite{CD89} or \cite{KMP}).
\end{remark}

\begin{remark}
Going back to the results for $L(p,q)$-spaces obtained in Section \ref{s2}, we
note that, for every $p_{0}>1$, the set ${\mathcal{A}}(p_{0},p_{0}%
)=\{(p,q):\,p\geq q\geq p_{0}\}$ consists precisely of the pairs of parameters
$(p,q)$, for which the space $L(p,q)$ is exact $p$-convex (see \cite{CD89} or
\cite{KMP}). Hence, on the one hand, we see that Corollary \ref{cor1ab}, which
gives sufficient conditions under which a Nazarov-Podkorytov type result holds
for $L(p,q)$-spaces, can be established in the special case when $p_{0}=q_{0}$
also by applying Theorem \ref{cor1c}(ii). On the other hand, Theorem
\ref{prop of Dlta}, which asserts that these conditions are also necessary,
shows that the assumption about the exact $p_{0}$-convexity of $X$ cannot be
removed from Theorem \ref{cor1c}(ii), even in the case of $L(p,q)$-spaces.

Moreover, recall that $L(p,q)$ is an interpolation space between $L^{p_{0}}$
and $L^{\infty}$ whenever $p>p_{0}$, $1\leq q\leq\infty$, or $p=q=p_{0}$ (see
e.g. \cite[Theorem~5.2.1]{BL}). Therefore, if $p,q$ satisfy these conditions,
by Theorem \ref{cor1c}(i), there exists a rearrangement invariant norm
$\Vert\cdot\Vert_{L(p,q)}^{\prime}$ on ${L(p,q),}$ equivalent to the initial
norm $\Vert\cdot\Vert_{L(p,q)},$ such that from the assumptions $f\in
L^{p_{0}}$, ${\Vert f\Vert}_{p_{0}}\geq{\Vert g\Vert}_{p_{0}}$, and
\eqref{distrib}, it follows that $\Vert f\Vert_{L(p,q)}^{\prime}\geq\Vert
g\Vert_{L(p,q)}^{\prime}$.
\end{remark}

\section{An extension of Nazarov-Podkorytov approach to rearrangement
invariant spaces\label{s8}}

Let us say that a pair of measurable functions $(f,g)$ satisfies the $NP$
condition (in brief, $(f,g)\in NP$) whenever condition \eqref{distrib1} holds
for some $t_{0}>0.$ In Section \ref{s7}, we have shown that if $(f,g)\in NP$
and, moreover, the inequality $\Vert f\Vert_{p}\geq\Vert g\Vert_{p}$ holds,
then $\Vert f\Vert_{Y}\geq\Vert g\Vert_{Y}$ for every rearrangement invariant
space $Y$ equipped with a norm of the form
\[
\Vert f\Vert_{Y}=\biggl\|\Bigl(\int_{0}^{t^{p}}f^{\ast}(s)^{p}\,ds\Bigr)^{1/p}%
\biggr\|_{F},
\]
where $F$ is a Banach function lattice on $(0,\infty)$. Furthermore, taking
into account equivalence \eqref{K-functLp} and the fact that $(L^{p}%
,L^{\infty})$ is a Calder\'{o}n-Mityagin couple, we deduced that this class of
spaces includes (up to equivalence of norms) all interpolation spaces with
respect to the pair $(L^{p},L^{\infty})$. It is natural to ask: what happens
when we replace $L^{p}$ by a rearrangement invariant space $X?$

We shall follow the analysis given in Section \ref{s7}. Recall that the
$K-$functional for the pair $(X,L^{\infty})$ satisfies, modulo absolute
constants of equivalence,
\[
K(t,f;X,L^{\infty})\asymp{\Vert f^{\ast}(\cdot)\chi_{(0,\phi_{X}^{-1}%
(t))}(\cdot)\Vert}_{X},\;\;t>0,
\]
where $\phi_{X}$ is the fundamental function of $X$, i.e., $\phi_{X}%
(t):=\Vert\chi_{(0,t)}(\cdot)\Vert_{X}$, and $\phi_{X}^{-1}(t)$ is its inverse
function\footnote{Without loss of generality, we can assume that $\phi_{X}$ is
a strictly increasing function on $(0,\infty)$.}. Let $\mathcal{X}$ be the
class of all rearrangement invariant spaces $X$ such that for every $(f,g)\in
NP$, with $\Vert f\Vert_{X}\geq\Vert g\Vert_{X},$ we can assert that%
\begin{equation}
{\Vert f^{\ast}(\cdot)\chi_{(0,s)}(\cdot)\Vert}_{X}\geq{\Vert g^{\ast}%
(\cdot)\chi_{(0,s)}(\cdot)\Vert}_{X},\text{ for all }s>0.
\label{extra condition}%
\end{equation}
If $X\in\mathcal{X}$, then for all $(f,g)\in NP$ such that $\Vert f\Vert
_{X}\geq\Vert g\Vert_{X},$ we immediately see that%
\[
\Vert f\Vert_{Y}\geq\Vert g\Vert_{Y}%
\]
for any rearrangement invariant space $Y$ defined by%
\[
\Vert h\Vert_{Y}=\biggl\|\Bigl\|h^{\ast}(\cdot)\chi_{(0,\phi_{X}^{-1}%
(t))}(\cdot)\Bigr\|_{X}\biggr\|_{F}<\infty,
\]
where $F$ is a Banach function lattice on $(0,\infty)$.

A complete characterization of the class $\mathcal{X}$ is an open problem.
Here we shall content ourselves with presenting examples of familiar spaces
that belong to the class $\mathcal{X}$, and others that do not.

We start with a positive result. Let $\varphi$ be an increasing concave
function on $[0,\infty)$ such that $\lim_{t\rightarrow0}\varphi(t)=0$ and let
$1\leq r<\infty$. The Lorentz space $\Lambda_{r}(\varphi)$ consists of all
measurable functions $f$ on $(0,\infty)$ such that
\[
{\Vert f\Vert}_{\Lambda_{r}(\varphi)}:=\Bigl(\int_{0}^{\infty}f^{\ast}%
(t)^{r}\,d\varphi(t)\Bigr)^{1/r}<\infty,
\]
(cf. \cite{Lo-51}, \cite[p. 121]{LT}, \cite{KM-04}).

\begin{proposition}
\label{Th5} Every Lorentz space $\Lambda_{r}(\varphi)$ belongs to the class
$\mathcal{X}$. In particular, $L(p,q)\in\mathcal{X}$ if $1<p<\infty$, $1\le
q\le p$.
\end{proposition}

\begin{proof}
To the contrary, let us suppose that for some $\varphi$ and $1\leq
r<\infty,\Lambda_{r}(\varphi)\not \in \mathcal{X}$. Then, there exists a pair
$(f,g)\in NP$ with $\Vert f\Vert_{\Lambda_{r}(\varphi)}\geq\Vert
g\Vert_{\Lambda_{r}(\varphi)},$ such that \eqref{extra condition} does not
hold: i.e., there is $s_{0}>0$ such that
\begin{equation}
{\Vert f^{\ast}(\cdot)\chi_{(0,s_{0})}(\cdot)\Vert}_{\Lambda_{r}(\varphi
)}<{\Vert g^{\ast}(\cdot)\chi_{(0,s_{0})}(\cdot)\Vert}_{\Lambda_{r}(\varphi)}.
\label{contrad1}%
\end{equation}
Let $t_{0}$ be the positive number associated with $(f,g)\in NP,$ then we must
have
\[
\Vert f^{\ast}(\cdot)\chi_{(0,s)}(\cdot)\Vert_{\Lambda_{r}(\varphi)}\geq\Vert
g^{\ast}(\cdot)\chi_{(0,s)}(\cdot)\Vert_{\Lambda_{r}(\varphi)}%
\]
if $s\leq t_{0}$. Hence, $s_{0}>t_{0}$. Therefore, by \eqref{distrib1}, for
$t>s_{0}$ we have $f^{\ast}(t)\leq g^{\ast}(t)$, which implies that
\[
\int_{s_{0}}^{\infty}f^{\ast}(t)^{r}\,d\varphi(t)\leq\int_{s_{0}}^{\infty
}g^{\ast}(t)^{r}\,d\varphi(t).
\]
Thus, by \eqref{contrad1},
\begin{align*}
{\Vert f\Vert}_{\Lambda_{r}(\varphi)}^{r}  &  =\int_{0}^{s_{0}}(f^{\ast
}(t))^{r}\,d\varphi(t)+\int_{s_{0}}^{\infty}(f^{\ast}(t))^{r}\,d\varphi(t)\\
&  ={\Vert f^{\ast}(\cdot)\chi_{(0,s_{0})}(\cdot)\Vert}_{\Lambda_{r}(\varphi
)}^{r}+\int_{s_{0}}^{\infty}(f^{\ast}(t))^{r}\,d\varphi(t)\\
&  <{\Vert g^{\ast}(\cdot)\chi_{(0,s_{0})}(\cdot)\Vert}_{\Lambda_{r}(\varphi
)}^{r}+\int_{s_{0}}^{\infty}(g^{\ast}(t))^{r}\,d\varphi(t)={\Vert g\Vert
}_{\Lambda_{r}(\varphi)}^{r},
\end{align*}
a contradiction. Consequently $\Lambda_{r}(\varphi)\in\mathcal{X}$. The second
part of the Proposition follows from the first since, for $1<p<\infty$, $1\leq
q\leq p$, we have, with equality of norms, $L(p,q)=\Lambda_{q}(\varphi)$,
where $\varphi(t)=t^{q/p}$.
\end{proof}

Let $\varphi$ be an increasing concave function on $[0,\infty)$ such that
$\lim_{t\rightarrow0}\varphi(t)=0$ and $\lim_{t\rightarrow\infty}%
\varphi(t)=\infty$. It is well known that (see e.g. \cite[Theorem~II.5.2]%
{KPS}) $\Lambda_{1}(\varphi)^{\ast}=M(\varphi)$, where $M(\varphi)$ is the
Marcinkiewicz space equipped with the norm
\[
{\Vert f\Vert}_{M(\varphi)}:=\sup_{t>0}\frac{1}{\varphi(t)}\int_{0}^{t}%
f^{\ast}(u)\,du.
\]
In contrast to Proposition \ref{Th5}, we have the following negative result
for Marcinkiewicz spaces.

\begin{proposition}
\label{Th5Marc} $M(t^{\alpha})\not \in \mathcal{X}$ for every $0<\alpha<1$.
\end{proposition}

We will need the following elementary lemma.

\begin{lemma}
\label{elementary} For each $0<\alpha<1$ it holds
\[
1-2^{-\alpha}<\alpha< 2(1-2^{-\alpha}).
\]

\end{lemma}

\begin{proof}
For the left-hand side inequality, it suffices to observe that the function
$F(\alpha):=\alpha-1+2^{-\alpha}$ increases for $\alpha\ge0$ and $F(0)=0$.

To check the remaining inequality, introduce the function $G(\alpha
):=2(1-2^{-\alpha})-\alpha$. Using calculus one can easily verify that $G$ has
a strict maximum at a point $\alpha_{0}\in(0,1)$: it increases on
$(0,\alpha_{0})$ and decreases on $(\alpha_{0},1)$. Since $G(0)=G(1)=0,$ the
desired result follows.
\end{proof}

\begin{proof}
[Proof of Proposition \ref{Th5Marc}]Let $0<\alpha<1,$ $\varphi(t)=t^{\alpha}.$
Consider the pair of functions $(f,g)$ defined by
\[
f(t)=\chi_{(0,1)}(t)+t^{-\alpha}\chi_{\lbrack1,\infty)}%
(t)\;\;\mbox{and}\;\;g(t)=\chi_{(0,2)}(t)+t^{-\alpha}\chi_{\lbrack2,\infty
)}(t),\;\;t>0.
\]
It is readily verified that the pair $(f,g)$ satisfies the $NP$ condition with
$t_{0}=1.$ To complete the proof we will show that ${\Vert f\Vert}%
_{M(\varphi)}={\Vert g\Vert}_{M(\varphi)}$, but there exists $s>0$ such that
inequality \eqref{extra condition} fails.

Observe that for $0<t<1$,
\[
\frac{\varphi(t)}{t}\int_{0}^{t}f(u)\,du=t^{\alpha},
\]
while for $t>1$ we have
\begin{equation}
\frac{\varphi(t)}{t}\int_{0}^{t}f(u)\,du=t^{\alpha-1}\left(  1+\frac
{1}{1-\alpha}(t^{1-\alpha}-1)\right)  =\frac{1}{1-\alpha}-\frac{\alpha
t^{\alpha-1}}{1-\alpha}. \label{wesee}%
\end{equation}
Hence,
\[
{\Vert f\Vert}_{M(\varphi)}=\lim_{t\rightarrow\infty}\frac{\varphi(t)}{t}%
\int_{0}^{t}f(u)\,du=\frac{1}{1-\alpha}.
\]
Also from (\ref{wesee}) we find
\[
{\Vert f\chi_{(0,2)}\Vert}_{M(\varphi)}=\frac{1}{1-\alpha}\left(
1-\alpha2^{\alpha-1}\right)  .
\]

Similarly, for $0<t<2$
\[
\frac{\varphi(t)}{t}\int_{0}^{t}g(u)\,du=t^{\alpha},
\]
and, for $t>2$,
\[
\frac{\varphi(t)}{t}\int_{0}^{t}g(u)\,du=t^{\alpha-1}\left(  2+\frac
{1}{1-\alpha}(t^{1-\alpha}-2^{1-\alpha})\right)  =\frac{1}{1-\alpha
}-2t^{\alpha-1}\left(  \frac{2^{-\alpha}}{1-\alpha}-1\right)  .
\]
By Lemma \ref{elementary}, $\frac{1}{1-\alpha}>2^{\alpha}$ or equivalently
$\frac{2^{-\alpha}}{1-\alpha}-1>0$. Therefore,
\[
{\Vert g\Vert}_{M(\varphi)}=\lim_{t\rightarrow\infty}\frac{\varphi(t)}{t}%
\int_{0}^{t}g(u)\,du=\frac{1}{1-\alpha}.
\]
Also, applying Lemma \ref{elementary} once again, we infer
\[
{\Vert g\chi_{(0,2)}\Vert}_{M(\varphi)}=2^{\alpha}>\frac{1}{1-\alpha}\left(
1-\alpha2^{\alpha-1}\right)  .
\]
Thus, ${\Vert f\Vert}_{M(\varphi)}={\Vert g\Vert}_{M(\varphi)}$ and
\[
{\Vert f\chi_{(0,2)}\Vert}_{M(\varphi)}<{\Vert g\chi_{(0,2)}\Vert}%
_{M(\varphi)}.
\]
Consequently, $M(\varphi)\not \in \mathcal{X}$.
\end{proof}

Our next result deals with Orlicz spaces.

An Orlicz function $M$ on $[0,\infty)$ satisfies the $\Delta_{2}$-condition
(briefly, $M\in\Delta_{2}$) whenever there exists a constant $C>0$ such that
\[
M(2u)\leq CM(u)\;\;\mbox{for all}\;\;u>0.
\]
Let $L_{M}$ be the Orlicz space $L_{M}$ equipped with the Luxemburg norm
defined by \eqref{luxemburg}.

\begin{theorem}
\label{Th6} Let $M$ be an Orlicz function on $[0,\infty)$ such that
$M\in\Delta_{2}$. Suppose that the function
\begin{equation}
\frac{M^{\prime}(\varepsilon x)}{M^{\prime}(x)}\quad
\mbox{ does not increase in $x$ on }(0,\infty)\;\mbox{for each }\varepsilon
\in(0,1). \label{condM}%
\end{equation}
Then, the Orlicz space $L_{M}$ belongs to the class $\mathcal{X}$.
\end{theorem}

\begin{proof}
Let $(f,g)\in NP$ for some $t_{0}>0$ and, moreover, $\Vert f\Vert_{L_{M}}%
\geq\Vert g\Vert_{L_{M}}$. Without loss of generality, we may assume that $f$
and $g$ are continuous on $(0,\infty)$ and $\Vert f\Vert_{L_{M}}=1$. Then,
since $M\in\Delta_{2}$, we have
\begin{equation}
\int_{0}^{\infty}M(g^{\ast}(t))\,dt\leq\int_{0}^{\infty}M(f^{\ast}%
(t))\,dt\leq1. \label{met1}%
\end{equation}

Suppose that $L_{M}\not \in \mathcal{X}$, then (cf. the proof of Proposition
\ref{Th5}) we have
\[
\label{contradiction1}{\Vert f^{\ast}(\cdot)\chi_{(0,s_{0})}(\cdot)\Vert
}_{L_{M}}<{\Vert g^{\ast}(\cdot)\chi_{(0,s_{0})}(\cdot)\Vert}_{L_{M}}%
\]
for some $s_{0}>t_{0}$. Therefore, using once again that $M\in\Delta_{2}$, we
can find $a>1$ such that
\begin{equation}
\int_{0}^{s_{0}}M(ag^{\ast}(t))\,dt>1\geq\int_{0}^{s_{0}}M(af^{\ast}(t))\,dt.
\label{met2}%
\end{equation}
Moreover, the same condition ensures that
\begin{equation}
\int_{0}^{s_{0}}M(ag^{\ast}(t))\,dt<\infty. \label{met2a}%
\end{equation}
Thus, to prove the desired result, it suffices to show that the hypothesis
\eqref{condM} of the theorem implies that inequalities \eqref{distrib1},
\eqref{met1} and \eqref{met2} are inconsistent. To this end, we need only to
prove that from the inequality
\[
\int_{0}^{s_{0}}M(ag^{\ast}(t))\,dt>\int_{0}^{s_{0}}M(af^{\ast}(t))\,dt,
\]
it follows that
\[
\int_{0}^{s_{0}}M(g^{\ast}(t))\,dt>\int_{0}^{s_{0}}M(f^{\ast}(t))\,dt.
\]
Indeed, since $g^{\ast}(t)\geq f^{\ast}(t)$ for $t>s_{0}>t_{0}$ (see
\eqref{distrib1}), the last inequality implies that
\[
\int_{0}^{\infty}M(g^{\ast}(t))\,dt>\int_{0}^{\infty}M(f^{\ast}(t))\,dt,
\]
which contradicts inequality \eqref{met1}.

Observe that condition \eqref{distrib1} holds as well for any scalar multiples
of the functions $f^{\ast}$ and $g^{\ast}$. Therefore, denoting $af^{\ast}$
and $ag^{\ast}$ by $F$ and $G$ respectively, it suffices to check that from
the inequality
\begin{equation}
\int_{0}^{s_{0}}M(G(t))\,dt>\int_{0}^{s_{0}}M(F(t))\,dt \label{met3}%
\end{equation}
it follows that
\begin{equation}
\int_{0}^{s_{0}}M(\varepsilon G(t))\,dt>\int_{0}^{s_{0}}M(\varepsilon
F(t))\,dt \label{met4}%
\end{equation}
for any $\varepsilon\in(0,1)$.

Further, we may assume that $F(t)=G(t)=0$ if $t>s_{0}$ (because the values of
$F(t)$ and $G(t)$ for $t>s_{0}$ play no role in the inequalities \eqref{met3}
and \eqref{met4}). Hence,
\[
\lim_{y\rightarrow0}M(y)\lambda_{F}(y)=\lim_{y\rightarrow0}M(y)\lambda
_{G}(y)=0.
\]
Moreover, taking into account \eqref{met2} and \eqref{met2a}, we see that the
functions $M(F)$ and $M(G)$ belong to the space $L^{1}(0,\infty)$, whence
\[
\lim_{t\rightarrow0}\int_{0}^{t}M(F(s))\,ds=\lim_{t\rightarrow0}\int_{0}%
^{t}M(G(s))\,ds=0,
\]
Therefore, since $M(F)$ and $M(G)$ are decreasing,
\[
\lim_{t\rightarrow0}M(F(t))t=\lim_{t\rightarrow0}M(G(t))t=0,
\]
and, passing to their distribution functions $\lambda_{F}$ and $\lambda_{G}$,
we get
\[
\lim_{y\rightarrow\infty}M(y)\lambda_{F}(y)=\lim_{y\rightarrow\infty
}M(y)\lambda_{G}(y)=0.
\]
Thus, setting $\tau_{0}:=F(t_{0})=G(t_{0})$ (see \eqref{distrib1}) and
integrating by parts, we can rewrite \eqref{met3} as follows:
\begin{equation}
\int_{0}^{\tau_{0}}(\lambda_{G}(y)-\lambda_{F}(y))M^{\prime}(y)\,dy>\int%
_{\tau_{0}}^{\infty}(\lambda_{F}(y)-\lambda_{G}(y))M^{\prime}(y)\,dy
\label{met5}%
\end{equation}
(note that both integrands are nonnegative on their integration domains).
Since for all functions $h,$ we have $\lambda_{\varepsilon h}(y)=\lambda
_{h}(y/\varepsilon)$, the inequality \eqref{met4} can be rewritten similarly:
\begin{equation}
\int_{0}^{\tau_{0}}(\lambda_{G}(y)-\lambda_{F}(y))M^{\prime}(\varepsilon
y)\,dy>\int_{\tau_{0}}^{\infty}(\lambda_{F}(y)-\lambda_{G}(y))M^{\prime
}(\varepsilon y)\,dy. \label{met6}%
\end{equation}
Then, by using \eqref{condM} and \eqref{met5}, we get
\begin{align*}
\int_{0}^{\tau_{0}}(\lambda_{G}(y)-\lambda_{F}(y))M^{\prime}(\varepsilon
y)\,dy  &  =\int_{0}^{\tau_{0}}(\lambda_{G}(y)-\lambda_{F}(y))\frac{M^{\prime
}(\varepsilon y)}{M^{\prime}(y)}M^{\prime}(y)\,dy\\
&  \geq\frac{M^{\prime}(\varepsilon\tau_{0})}{M^{\prime}(\tau_{0})}\int%
_{0}^{\tau_{0}}(\lambda_{G}(y)-\lambda_{F}(y))M^{\prime}(y)\,dy\\
&  >\frac{M^{\prime}(\varepsilon\tau_{0})}{M^{\prime}(\tau_{0})}\int_{\tau
_{0}}^{\infty}(\lambda_{G}(y)-\lambda_{F}(y))M^{\prime}(y)\,dy\\
&  \geq\int_{\tau_{0}}^{\infty}(\lambda_{G}(y)-\lambda_{F}(y))\frac{M^{\prime
}(\varepsilon y)}{M^{\prime}(y)}M^{\prime}(y)\,dy\\
&  =\int_{\tau_{0}}^{\infty}(\lambda_{G}(y)-\lambda_{F}(y))M^{\prime
}(\varepsilon y)\,dy,
\end{align*}
and \eqref{met6} follows, completing the proof.
\end{proof}

The function $M(t)=t^{p}$, $1\leq p<\infty$, is the simplest example of an
Orlicz function that satisfies all conditions of Theorem \ref{Th6}. Indeed,
for each $\varepsilon>0$ we have
\[
\frac{M^{\prime}(\varepsilon x)}{M^{\prime}(x)}=\frac{p(\varepsilon x)^{p-1}%
}{px^{p-1}}=\varepsilon^{p-1},
\]
and hence \eqref{condM} holds.

Let us now present some non-trivial examples of Orlicz functions which satisfy
the conditions of Theorem \ref{Th6}.

\begin{example}
Let $0<p<q$ and let $N(u)=\max\{u^{p},u^{q}\}$ for $u\in\lbrack0,\infty)$.
Therefore $N(u)$ is positive and increasing on $(0,\infty)$, and
\[
M(x):=\int_{0}^{x}N(u)\,du=%
\begin{cases}
\quad\quad\;\;\frac{x^{p+1}}{p+1}\quad\quad\quad\mbox{if }x\in\lbrack0,1],\\
\frac{q-p}{(p+1)(q+1)}+\frac{x^{q+1}}{q+1}\;\;\mbox{if }x\in(1,\infty),
\end{cases}
\]
is an Orlicz function such that $M\in\Delta_{2}$. Moreover, for every
$\varepsilon\in(0,1)$
\[
\frac{M^{\prime}(\varepsilon x)}{M^{\prime}(x)}=\frac{N(\varepsilon x)}{N(x)}=%
\begin{cases}
\;\;\varepsilon^{p}\quad\quad\mbox{if }x\in(0,1),\\
\varepsilon^{p}x^{p-q}\;\;\mbox{if }x\in\lbrack1,1/\varepsilon],\\
\;\;\varepsilon^{q}\quad\quad\mbox{if }x\in(1/\varepsilon,\infty)
\end{cases}
\]
is non-increasing on $(0,\infty)$. Therefore, $M$ satisfies the conditions of
Theorem \ref{Th6}, and therefore, $L_{M}\in\mathcal{X}$.
\end{example}

\begin{example}
Let $N(u)$ be defined on $[0,\infty)$ by
\[
N(u)=%
\begin{cases}
\;\;u\quad\;\;\mbox{if }u\in\lbrack0,1],\\
\frac{u^{2}}{\ln(eu)}\;\;\mbox{if }u\in(1,\infty).
\end{cases}
\]
Then, $N(u)$ is positive, and since
\[
N^{\prime}(u)=%
\begin{cases}
\;\;\quad1\quad\quad\quad\mbox{if }u\in\lbrack0,1],\\
\frac{u(2\ln(eu)-1)}{\ln^{2}(eu)}\;\;\mbox{if }u\in(1,\infty),
\end{cases}
\]
we see that $N(u)$ is increasing on $(1,\infty)$. It follows that,
\[
M(x)=\int_{0}^{x}N(u)\,du
\]
is an Orlicz function. One can easily verify that $M(x)$ satisfies the
$\Delta_{2}$-condition. Indeed, for $x>1,$
\[
M(2x)=M(x)+\int_{x}^{2x}N(u)\,du\leqslant M(x)+\frac{1}{\ln(ex)}\int_{x}%
^{2x}u^{2}\,du=M(x)+\frac{7x^{3}}{3\ln(ex)},
\]
and
\[
M(x)\geqslant\int_{1}^{x}N(u)\,du\geqslant\frac{1}{\ln(ex)}\int_{1}^{x}%
u^{2}\,du=\frac{x^{3}-1}{3\ln(ex)}.
\]
Consequently, there exists $C>0,$ such that for $x$ sufficiently large we have
$M(2x)\leqslant CM(x)$. Since $M(x)=\frac{x^{2}}{2}$, $x\in\lbrack0,1]$, see
that there exists $C_{1}>0,$ such that $M(x)$ for some $C_{1}>0$ and all
$x\in(0,\infty)$.

Next, for each $\varepsilon\in(0,1)$ the function
\[
\frac{M^{\prime}(\varepsilon x)}{M^{\prime}(x)}=\frac{N(\varepsilon x)}{N(x)}=%
\begin{cases}
\;\;\quad\varepsilon\quad\quad\quad\mbox{if }x\in(0,1),\\
\quad\frac{\varepsilon\ln(ex)}{x}\;\;\quad\mbox{if }x\in\lbrack1,1/\varepsilon
],\\
\;\;\frac{\varepsilon^{2}\ln(ex)}{\ln(ex)+\ln\varepsilon}\quad\mbox{if }x\in
(1/\varepsilon,\infty)
\end{cases}
\]
is nonincreasing on $(0,\infty)$. Then, by Theorem \ref{Th6}, the Orlicz space
$L_{M}$ belongs to the class $\mathcal{X}$.
\end{example}

\begin{example}
Similarly, if we let
\[
N(u)=%
\begin{cases}
\frac{u}{\ln(e/u)}\quad\;\;\mbox{if }u\in\lbrack0,1],\\
\quad{u^{2}}\quad\quad\mbox{if }u\in(1,\infty),
\end{cases}
\]
then the Orlicz function $M,$ defined by
\[
M(x)=\int_{0}^{x}N(u)\,du,
\]
satisfies the conditions of Theorem \ref{Th6}, and therefore $L_{M}%
\in\mathcal{X}$.
\end{example}

\section{Comparison of norms: differential inequalities and
extrapolation\label{differential inequalities}}

\subsection{Introduction\label{Introd}}

As we have shown through this paper, the Nazarov-Podkorytov Lemma can be
naturally placed in the context of the more general problem of comparing norms
of elements that belong to an interpolation scale. Let $\{X_{\theta}%
\}_{\theta\in(0,1)}$ be an interpolation scale and let $f,g$ be elements, for
which we know, for example, that for some $\theta_{0}\in(0,1)$
\[
\left\Vert f\right\Vert _{X_{\theta_{0}}}=\left\Vert g\right\Vert
_{X_{\theta_{0}}}.
\]
What other conditions do we need to impose in order to ascertain that there
exists a constant $C$ (independent of $f$ and $g)$ such that for some
$\theta_{1}\neq\theta_{0}$
\[
\left\Vert g\right\Vert _{X_{\theta_{1}}}\leq C\left\Vert f\right\Vert
_{X_{\theta_{1}}}?
\]
As we have seen, in the case of real interpolation scales it is natural to
formulate conditions on the underlying $K-$ or $E-$functionals, which, in
turn, are generalizations of the classical real analysis concepts of
distribution functions and rearrangements, etc. Here, we indicate a different
set of conditions to compare the norms of $f$ and $g$ which is based on
comparing the derivatives of their norms as we move through the scale (i.e.,
as $\theta$ changes). In this fashion the comparison criteria that we
formulate will be based on the fundamental theorem of calculus!

Moreover, as we will show, there are mechanisms to use these ideas to obtain
estimates on the underlying $K-$ and $E-$functionals themselves. Needless to
say that to fully develop this topic falls outside the scope of the present
work, but we feel it is important to indicate some basic results and their
connection to the Nazarov-Podkorytov approach.

Indeed, taking derivatives of the norms appeared early on in the theory of
Banach spaces. Of particular interest to us are computations that involve a
scale of interpolation spaces. For example, in \cite{mazur} it was shown that
the $L^{p}$-norm is Frechet differentiable
and its directional derivatives were computed.

\begin{example}
Let $\Omega$ be a probability space and let $g\in L^{q}(\Omega).$ Then for
$p<q,$ we have (cf. \cite{mazur}, \cite[(2.6), page 1065]{gross})
\begin{equation}
\frac{d}{dp}(\left\Vert g\right\Vert _{p})=\frac{1}{p}\left\Vert g\right\Vert
_{p}^{1-p}\cdot\int_{\Omega}\left\vert g\right\vert ^{p}\ln\frac{\left\vert
g\right\vert }{\left\Vert g\right\Vert _{p}}\,ds. \label{entro}%
\end{equation}

\end{example}

In this section we consider a method to compare norms in interpolation scales
via differential inequalities. Related ideas but framed in a slightly
different language appeared in some form already in \cite{cjm} and \cite{cjmr}
but at that time it was not clear what was the direction to take concerning
applications and connections with other topics. We proceed to informally
explain the simple underlying ideas. Let $\{X_{\theta}\}_{\theta\in(0,1)}$ be
an ordered interpolation scale of Banach spaces in the sense that
$X_{1}\subset X_{\theta_{1}}\subset X_{\theta_{0}}\subset X_{0},$ when
$0<\theta_{0}<\theta_{1}<1$. In particular, under this assumption, we have
$\frac{d}{d\theta}\left\Vert \cdot\right\Vert _{X_{\theta}}\geq0,$ when it
makes sense. Suppose that $0<\theta_{0}<\theta_{1}<1$ are fixed, and $f,g$ are
elements that belong to $X_{\theta_{1}}$ such that%
\begin{equation}
\left\Vert g\right\Vert _{X_{\theta_{0}}}\leq\left\Vert f\right\Vert
_{X_{\theta_{0}}}. \label{inview}%
\end{equation}
Our aim, just like in other parts of this paper, is to estimate $\left\Vert
g\right\Vert _{X_{\theta_{1}}}$ in terms of $\left\Vert f\right\Vert
_{X_{\theta_{1}}}.$ Suppose that we can compare $\frac{d}{d\theta}(\left\Vert
g\right\Vert _{X_{\theta}})$ with $\frac{d}{d\theta}(\left\Vert f\right\Vert
_{X_{\theta}})$ as follows: There exists a constant $C>0$ such that
\[
\frac{d}{d\theta}(\left\Vert g\right\Vert _{X_{\theta}})\leq C\frac{d}%
{d\theta}(\left\Vert f\right\Vert _{X_{\theta}}),\text{ for all }\theta
_{0}<\theta<\theta_{1}.
\]
Then, using the Fundamental Theorem of Calculus, we find%
\begin{align*}
\left\Vert g\right\Vert _{X_{\theta_{1}}}-\left\Vert g\right\Vert
_{X_{\theta_{0}}}  &  =\int_{\theta_{0}}^{\theta_{1}}\frac{d}{d\theta
}(\left\Vert g\right\Vert _{X_{\theta}})\,ds\\
&  \leq C\int_{\theta_{0}}^{\theta_{1}}\frac{d}{d\theta}(\left\Vert
f\right\Vert _{X_{\theta}})\,ds\\
&  =C\left(  \left\Vert f\right\Vert _{X_{\theta_{1}}}-\left\Vert f\right\Vert
_{X_{\theta_{0}}}\right)  .
\end{align*}
Therefore, if $C\geq1,$ then in view of (\ref{inview}),
\begin{align*}
\left\Vert g\right\Vert _{X_{\theta_{1}}}  &  \leq C\left\Vert f\right\Vert
_{X_{\theta_{1}}}+\underset{\leq0}{\underbrace{\left(  \left\Vert g\right\Vert
_{X_{\theta_{0}}}-C\left\Vert f\right\Vert _{X_{\theta_{0}}}\right)  }}\\
&  \leq C\left\Vert f\right\Vert _{X_{\theta_{1}}}.
\end{align*}

One can conceive the process of integrating differential inequalities we have
described as an extrapolation result, where we obtain an end point result from
a family of estimates. Indeed, for example, one could replace the constant $C$
by a weight $C(\theta)$ and apply the mechanisms of extrapolation theory (cf.
\cite{AM}, \cite{ALM}) to obtain a more general family of such results.

In the next section we show a connection between norm comparisons, derivatives
of norms, $K-$functionals and certain operators associated to their computation.

\subsection{$\Omega$ operators and derivatives of norms}

\label{derivatives of norms}

Recall that we work with ordered pairs of Banach spaces $(X_{0},X_{1}),$ that
is, we assume that with continuous embedding,%
\begin{equation}
X_{1}\subset X_{0}. \label{no0}%
\end{equation}
In this setting sometimes it is convenient to renormalize the definition of
the interpolation spaces. For $\theta\in(0,1),1\leq q\leq\infty,$ we define
(cf. \cite{ALM} and references therein)%
\begin{equation}
\lbrack X_{0},X_{1}]_{\theta,q}=\{f\in X_{0}:\left\Vert f\right\Vert
_{[X_{0},X_{1}]_{\theta,q}}=\left\{  \int_{0}^{1/2}[s^{-\theta}K(s,f;X_{0}%
,X_{1})]^{q}\frac{ds}{s}\right\}  ^{1/q}<\infty\}, \label{no1}%
\end{equation}
with the usual modification when $q=\infty.$ The difference between the usual
norm $\left\Vert f\right\Vert _{(X_{0},X_{1})_{\theta,q}}$ and $\left\Vert
f\right\Vert _{[X_{0},X_{1}]_{\theta,q}}$ is \ that in (\ref{no1}) the
interval of integration is $(0,1/2)$ rather than $(0,\infty).$ However, since
we are under the regime (\ref{no0}) the $K-$functional becomes constant for
$t$ large enough (say $t$ bigger than the norm of the embedding (\ref{no0}))
and therefore the classical and the modified norms\footnote{On the other hand,
the modified spaces $[X_{0},X_{1}]_{\theta,q}$ are useful to define limiting
cases, e.g. $\theta=0,$ for which the usual norms can only be finite on the
null element.} are equivalent.

Associated with the computation of the $K-$, $E-$, $J-$functionals are certain
operators that also play a r\^{o}le in the calculations below. We present a
brief summary of the relevant definitions and results for the convenience of
the reader.

Suppose that $f\in X_{0},$ and for each $t\in(0,1)$ select a decomposition%
\[
f=a_{0}(t)+a_{1}(t),\;a_{i}(t)\in X_{i},\;\;i=0,1,
\]
such that with constants of equivalence independent of $t,$%
\[
K(t,f;X_{0},X_{1})\asymp\left\Vert a_{0}(t)\right\Vert _{X_{0}}+t\left\Vert
a_{1}(t)\right\Vert _{X_{1}},\;\;t\in(0,1).
\]
The $\Omega$ operator is defined (modulo bounded operators) by%
\[
\Omega f:=\Omega_{(X_{0},X_{1})}f=\int_{0}^{1/2}a_{0}(s)\frac{ds}{s}%
-\int_{1/2}^{1}a_{1}(s)\frac{ds}{s}.
\]
Then, for any $t\in(0,1/2)$%
\begin{align*}
\Omega f  &  =\int_{0}^{t}a_{0}(s)\frac{ds}{s}+\int_{t}^{1/2}a_{0}(s)\frac
{ds}{s}-\int_{1/2}^{1}a_{1}(s)\frac{ds}{s}\\
&  =\int_{0}^{t}a_{0}(s)\frac{ds}{s}+\int_{t}^{1/2}(a_{0}(s)+a_{1}%
(s))\frac{ds}{s}-\int_{t}^{1}a_{1}(s)\frac{ds}{s}\\
&  =\int_{0}^{t}a_{0}(s)\frac{ds}{s}-\int_{t}^{1}a_{1}(s)\frac{ds}{s}-f\log2t.
\end{align*}

It follows that\footnote{Since by definition $\left\Vert a_{0}(s)\right\Vert
_{X_{0}}\leq K(s,f),$and $\left\Vert a_{1}(s)\right\Vert _{X_{1}}=\frac{s}%
{s}\left\Vert a_{1}(s)\right\Vert _{X_{1}}\leq\frac{K(s,f)}{s}.$}%
\begin{align*}
K(t,\Omega f+f\log2t) &  \leq\left\Vert \int_{0}^{t}a_{0}(s)\frac{ds}%
{s}\right\Vert _{X_{0}}+t\left\Vert \int_{t}^{1}a_{1}(s)\frac{ds}%
{s}\right\Vert _{X_{1}}\\
&  \leq\int_{0}^{t}\left\Vert a_{0}(s)\right\Vert _{X_{0}}\frac{ds}{s}%
+t\int_{t}^{1}\left\Vert a_{1}(s)\right\Vert _{X_{1}}\frac{ds}{s}\\
&  \preceq\int_{0}^{1}K(s,f)\min\{1,\frac{t}{s}\}\frac{ds}{s}.
\end{align*}
In particular, if $T$ is a bounded operator, $T:X_{i}\rightarrow Y_{i},i=0,1,$
then letting $\left[  \Omega,T\right]  f:=T\Omega_{(X_{0},X_{1})}%
f-\Omega_{(Y_{0},Y_{1})}Tf$ we have
\begin{align*}
K(t,\left[  \Omega,T\right]  f;Y_{0},Y_{1}) &  \leq K(t,T(\Omega_{(X_{0}%
,X_{1})}f-f\log2t);Y_{0},Y_{1})\\
&  +K(t,(\log2t)Tf-\Omega_{(Y_{0},Y_{1})}Tf;Y_{0},Y_{1})\\
&  \preceq\left\Vert T\right\Vert \int_{0}^{1}K(s,f;X_{0},X_{1})\min
\{1,\frac{t}{s}\}\frac{ds}{s}\\
&  \preceq\left\Vert T\right\Vert S(K(\cdot,f;X_{0},X_{1}))(t),
\end{align*}
where $S$ is the Calder\'{o}n operator defined by $S\phi(t)=\int_{0}^{1}%
\phi(s)\min\{1,\frac{t}{s}\}\frac{ds}{s}.$ Consequently, if we let
\[
\left\Vert \phi\right\Vert _{L_{\theta}^{q}}=\left\{  \int_{0}^{1}%
[\phi(t)t^{-\theta}]^{q}\frac{dt}{t}\right\}  ^{1/q}<\infty,\;\;0<\theta
<1,\;1\leq q<\infty,
\]
we can select $c(\theta,q)$ such that for all $\phi,$
\[
\left\Vert S\phi\right\Vert _{L_{\theta}^{q}}\leq c(\theta,q)\left\Vert
\phi\right\Vert _{L_{\theta}^{q}},
\]
and we obtain
\begin{equation}
\left\Vert \left[  \Omega,T\right]  f\right\Vert _{[Y_{0},Y_{1}]_{\theta,q}%
}\preceq\left\Vert T\right\Vert c(\theta,q)\left\Vert f\right\Vert
_{[X_{0},X_{1}]_{\theta,q}}.\label{normd}%
\end{equation}
For some calculations it is preferable to use $E-$, $J-$functionals or some
variants of these and $K-$functionals, for which corresponding versions of
(\ref{normd}) also hold (cf. \cite{jrw}, \cite{cjmr}). For example, for the
$E-$method, if $0<\alpha<\infty$, $q>0$ and%
\[
(X_{0},X_{1})_{\alpha,q;E}=\{f:\left\Vert f\right\Vert _{(X_{0},X_{1}%
)_{\theta,q;E}}=\left\{  \int_{0}^{\infty}[s^{\alpha}E(s,f,X_{0},X_{1}%
)]^{q}\frac{ds}{s}\right\}  ^{1/q}<\infty,
\]
for $\theta=1/(\alpha+1)$ we have (cf. \cite[Theorem 7.1.7, page 178]{BL},
\cite{cjmr})%
\[
\left\Vert f\right\Vert _{(X_{0},X_{1})_{\frac{1-\theta}{\theta},\theta q;E}%
}=\left\Vert f\right\Vert _{[X_{0},X_{1}]_{\theta,q}}^{1/\theta}.
\]

To fix ideas we only consider the case $q=1$ and let
\[
\left\Vert f\right\Vert _{X_{\theta}}=\left\Vert f\right\Vert _{[X_{0}%
,X_{1}]_{\theta,1}}.
\]
Then (cf. \cite{cjm}, \cite{cjmr})%
\begin{align*}
\frac{d}{d\theta}(\left\Vert h\right\Vert _{X_{\theta}})  &  =\int_{0}%
^{1/2}[s^{-\theta}\log\frac{1}{s}K(s,h;X_{0},X_{1})]\frac{ds}{s}\\
&  \asymp\left\Vert h\right\Vert _{X_{\theta}}+\left\Vert \Omega h\right\Vert
_{X_{\theta}}.
\end{align*}

In particular, conditions of the form
\begin{equation}
\frac{d}{d\theta}(\left\Vert g\right\Vert _{X_{\theta}})\leq C\frac{d}%
{d\theta}(\left\Vert f\right\Vert _{X_{\theta}})\text{ uniformly for all
}\theta\in(\theta_{0},\theta_{1}) \label{deriva}%
\end{equation}
are equivalent to
\begin{align*}
\int_{0}^{1/2}[s^{-\theta}\log\frac{1}{s}K(s,g;X_{0},X_{1})]\frac{ds}{s}  &
\leq\int_{0}^{1/2}[s^{-\theta}\log\frac{1}{s}K(s,f;X_{0},X_{1})]\frac{ds}%
{s},\\
\text{ uniformly for all }\theta &  \in(\theta_{0},\theta_{1}).
\end{align*}
Consequently, we have a different way to formulate condition (\ref{deriva}),
which combined with an inequality of the form $\left\Vert g\right\Vert
_{X_{\theta_{0}}}\leq C\left\Vert f\right\Vert _{X_{\theta_{0}}},$ allows to
conclude that
\[
\left\Vert g\right\Vert _{X_{\theta_{1}}}\leq C\left\Vert f\right\Vert
_{X_{\theta_{1}}}.
\]
By choosing suitable $E-$ or $J-$methods we can easily compute derivatives of
more \textquotedblleft complicated\textquotedblright\ norms, but again here we
just wanted to present the basic ideas of the method.

\section{Appendix}

\subsection{A simplified version of the Nazarov-Podkorytov Lemma}

\label{s9}

We prove an elementary result showing that inequalities of the type $\Vert
f\Vert_{p}\geq\Vert g\Vert_{p}$, for $p$ sufficiently large, can be obtained
under essentially weaker conditions than those imposed in the
Nazarov-Podkorytov Lemma.

\begin{proposition}
\label{cor1d} Let $1\leq p_{0}<\infty$, $0<t_{0}<\infty$, and let $f,g$ be
such that $f,g\in L^{p_{0}}$, $f^{\ast}(t_{0})>g^{\ast}(t_{0})$ and $f^{\ast
}(t)\geq g^{\ast}(t)$ for all $0<t<t_{0}$. Suppose that the set $\mathcal{P}%
:=\{p>0:\,|f|^{p}-|g|^{p}\in L^{1}(0,\infty)\}$ is unbounded. Then, for all
$p\in\mathcal{P}$, that are large enough, we have
\[
\int_{0}^{\infty}(|f(t)|^{p}-|g(t)|^{p})\,dt>0.
\]

\end{proposition}

\begin{proof}
Suppose that $h\in L^{p_{0}}\cap L^{\infty}$. Then, $h\in L^{p}$ for all
$p\in\lbrack p_{0},\infty)$. Let us prove that%
\begin{equation}
\lim_{p\rightarrow\infty}{\Vert h\Vert}_{p}={\Vert h\Vert}_{\infty}.
\label{limit}%
\end{equation}
From the fact that $h\in L^{p_{0}}$, it follows that for any $\varepsilon>0$
there exists $t_{1}=t_{1}(\varepsilon)$ such that $h^{\ast}(t)<\varepsilon$
for all $t>t_{1}$. Consequently, for $p>p_{0}$,
\begin{align*}
{\Vert h\Vert}_{p}  &  \leq\biggl(\int_{0}^{t_{1}}(h^{\ast}(s))^{p}%
\,ds\biggr)^{1/p}+\biggl(\int_{t_{1}}^{\infty}(h^{\ast}(s))^{p}%
\,ds\biggr)^{1/p}\\
&  \leq{\Vert h\Vert}_{\infty}\cdot\left(  t_{1}^{1/p}+\biggl(\int_{t_{1}%
}^{\infty}(h^{\ast}(s))^{p_{0}}\cdot\varepsilon^{p-p_{0}}\,ds\biggr)^{1/p}%
\right) \\
&  \leq{\Vert h\Vert}_{\infty}\cdot\left(  t_{1}^{1/p}+\varepsilon
^{\frac{p-p_{0}}{p}}{\Vert h\Vert}_{p_{0}}^{p_{0}/p}\right)  .
\end{align*}
Since the right-hand side of this inequality tends to ${\Vert h\Vert}_{\infty
}(1+\varepsilon)$ as $p\rightarrow\infty$, we get
\[
\underset{p\rightarrow\infty}{\lim\sup}{\Vert h\Vert}_{p}\leq{\Vert h\Vert
}_{\infty}.
\]

From the fact that $h\in L^{\infty},$ we can find $t_{2}>0$, such that
$h^{\ast}(t_{2})>{\Vert h\Vert}_{\infty}-\varepsilon$. Hence,
\[
{\Vert h\Vert}_{p}\geq\biggl(\int_{0}^{t_{2}}(h^{\ast}(s))^{p}%
\,ds\biggr)^{1/p}\geq({\Vert h\Vert}_{\infty}-\varepsilon)t_{2}^{1/p}.
\]
Therefore, passing to the limit as $p\rightarrow\infty$, we arrive at the
inequality
\[
{\Vert h\Vert}_{p}\geq{\Vert h\Vert}_{\infty}-\varepsilon.
\]
Thus,
\[
\underset{p\rightarrow\infty}{\lim\inf}{\Vert h\Vert}_{p}\geq{\Vert h\Vert
}_{\infty}%
\]
Combining these inequalities, we establish \eqref{limit}.

Let $p\in\mathcal{P}$. Combining in a familiar manner the fact that
$|f|^{p}-|g|^{p}\in L^{1}(0,\infty)$, and Fubini's theorem (cf. the proof of
Proposition \ref{cor1b}), yields
\begin{equation}
\int_{0}^{\infty}(|f(t)|^{p}-|g(t)|^{p})\,dt=\int_{0}^{\infty}(f^{\ast}%
(t)^{p}-g^{\ast}(t)^{p})\,dt. \label{EQ1-App}%
\end{equation}

By \cite[Theorem~II.3.1]{KPS} it follows that $f^{\ast}(t)^{p}-g^{\ast}%
(t)^{p}$ belongs to $L^{1}(0,\infty)$ together with $|f(t)|^{p}-|g(t)|^{p}$,
and by assumption,
\begin{equation}
\int_{0}^{t_{0}}(f^{\ast}(s)^{p}-g^{\ast}(s)^{p})\,ds\geq0. \label{EQ2-App}%
\end{equation}

Pick $\varepsilon>0$ so that $f^{\ast}(t_{0})>g^{\ast}(t_{0})+2\varepsilon$,
then applying \eqref{limit} to the functions $f^{\ast}\cdot\chi_{(t_{0}%
,\infty)}$ and $g^{\ast}\cdot\chi_{(t_{0},\infty)}$, we see that for all $p$
large enough
\[
{\Vert f^{\ast}\cdot\chi_{(t_{0},\infty)}\Vert}_{p}>f^{\ast}(t_{0}%
)-\varepsilon>g^{\ast}(t_{0})+\varepsilon>{\Vert g^{\ast}\cdot\chi
_{(t_{0},\infty)}\Vert}_{p}.
\]
Therefore, by \eqref{EQ1-App} and \eqref{EQ2-App}, for those values of
$p\in\mathcal{P}$ we have
\begin{align*}
\int_{0}^{\infty}(|f(t)|^{p}-|g(t)|^{p})\,dt  &  =\int_{0}^{t_{0}}(f^{\ast
}(t)^{p}-g^{\ast}(t)^{p})\,dt+\int_{t_{0}}^{\infty}(f^{\ast}(t)^{p}-g^{\ast
}(t)^{p})\,dt\\
&  \geq{\Vert f^{\ast}\cdot\chi_{(t_{0},\infty)}\Vert}_{p}^{p}-{\Vert g^{\ast
}\cdot\chi_{(t_{0},\infty)}\Vert}_{p}^{p}>0.
\end{align*}

\end{proof}

\begin{remark}
The conditions of Proposition \ref{cor1d} are weaker than those of the
Nazarov-Podkorytov Lemma but do not allow to determine a value $p_{0}$ for
which ${\Vert f\Vert}_{p}\geq{\Vert g\Vert}_{p}$ holds for $p\geq p_{0}$.
\end{remark}

\subsection{Modular inequalities of Nazarov-Podkorytov type and Karamata
inequality.\label{s10}}


To conclude we show that there are close connections between the
Nazarov-Podkorytov Lemma and the well-known classical Karamata inequality.
More specifically, using the latter inequality, one can easily obtain rather
general Nazarov-Podkorytov type results for modulars.

We recall a discrete version of Karamata's inequality (see e.g. \cite{Kar},
\cite{HLP}, \cite[Ch.~I, \S \,28]{BB}).
%
Suppose that ${x}=(x_{n})_{n=1}^{\infty}$, ${y}=(y_{n})_{n=1}^{\infty}$ are
two nonincreasing sequences of real numbers, and let $\Phi=\Phi(u)$ be an
arbitrary nondecreasing convex function. Then, from
\begin{equation}
\sum_{j=1}^{k}x_{j}\geqslant\sum_{j=1}^{k}y_{j},\quad
\mbox{ for all }k=1,2,\ldots\label{SeqIneqKaramataInf}%
\end{equation}
it follows that%
\[
\sum_{j=1}^{\infty}\Phi(x_{j})\geqslant\sum_{j=1}^{\infty}\Phi(y_{j}).
\]

It is now easy to see that the usual cancellation argument implies a rather
general discrete Nazarov-Podkorytov type result for modulars. Recall (see
Section \ref{s4}), that if ${x}=(x_{n})_{n=1}^{\infty}$ and ${y}=(y_{n}%
)_{n=1}^{\infty}$ are two sequences, we write ${x}\gtrless{y}$, whenever for
some $n_{0}\in\mathbb{N}$
\[
x_{n}^{\ast}\geqslant y_{n}^{\ast}\;\;\mbox{if}\;\;n\leqslant n_{0}%
\;\;\mbox{and}\;\;x_{n}^{\ast}\leqslant y_{n}^{\ast}\;\;\mbox{if}\;\;n>n_{0}.
\]

\begin{proposition}
Suppose that $\Phi$ is a nondecreasing function and $\Psi$ is a nondecreasing
and convex function on $[0,\infty)$. Then, for any sequences ${x}%
=(x_{n})_{n=1}^{\infty}$ and ${y}=(y_{n})_{n=1}^{\infty}$, which satisfy the
conditions: ${x}\gtrless{y}$ and
\[
\sum_{j=1}^{\infty}\Phi(|x_{j}|)\geqslant\sum_{j=1}^{\infty}\Phi(|y_{j}|),
\]
we have
\[
\sum_{j=1}^{\infty}\Psi(\Phi(|x_{j}|))\geqslant\sum_{j=1}^{\infty}\Psi
(\Phi(|y_{j}|)).
\]

\end{proposition}

\begin{proof}
Since $\Phi$ is nondecreasing, the nonincreasing permutation of the sequence
$(\Phi(|x_{n}|))_{n=1}^{\infty}$ (resp. $(\Phi(|y_{n}|))_{n=1}^{\infty}$)
coincides with $(\Phi(x_{n}^{\ast}))_{n=1}^{\infty}$ (resp. $(\Phi(y_{n}%
^{\ast}))_{n=1}^{\infty}$). Proceeding as in Section \ref{s7} (see Lemma
\ref{lemma:informal}), it follows that these sequences satisfy a condition
analogous to \eqref{SeqIneqKaramataInf}, i.e.,
\[
\sum_{n=1}^{k}\Phi(x_{n}^{\ast})\geqslant\sum_{n=1}^{k}\Phi(y_{n}^{\ast}%
)\quad\mbox{ for all }k=1,2,\ldots.
\]
Therefore, applying Karamata's inequality to the sequences $(\Phi(x_{n}^{\ast
}))_{n=1}^{\infty}$, $(\Phi(y_{n}^{\ast}))_{n=1}^{\infty}$ and the function
$\Psi$, we obtain the desired inequality.
\end{proof}

Karamata's inequality can be also formulated for measurable functions defined
on $(0,\infty)$ (see e.g. \cite{HLP}, \cite{KFL}, \cite[Ch.~I, \S \,32]{BB}).
As a consequence we obtain

\begin{proposition}
Suppose that $\Phi$ is a nondecreasing function and $\Psi$ is a nondecreasing
and convex function on $[0,\infty)$. Then, for any pair of functions $f$ and
$g$, which satisfy condition \eqref{distrib1}, for some $t_{0}>0,$ and
moreover, such that
\[
\int_{0}^{\infty}\Phi(|f(t)|)\,dt\geqslant\int_{0}^{\infty}\Phi(|g(t)|)\,dt,
\]
we have
\[
\int_{0}^{\infty}\Psi(\Phi(|f(t)|))\,dt\geqslant\int_{0}^{\infty}\Psi
(\Phi(|g(t)|))\,dt.
\]

\end{proposition}

A special case of the previous proposition, i.e., letting $\Phi(u)=u^{p}$,
$\Psi(u)=u^{q/p}$, $q>p>0$ gives the Nazarov-Podkorytov Lemma.

\end{document}